\documentclass[12pt]{article}
\usepackage{epsfig,amssymb,epic,eepic,color,graphics,amsmath}
\usepackage[utf8]{inputenc}
\usepackage{multicol, color}
\input{amssym.def}

\title{There are no structural stable Axiom A 3-diffeomorphisms with dynamics ``one-dimensional surfaced attractor-repeller''}
\author{O. Pochinka}
\date{DSA Laboratory HSE University (Russia)}

\usepackage{graphicx}
\newtheorem{theo}{Theorem}

\newtheorem{lemm}{Lemma}

\newenvironment{demo}{{\bf Proof: }}{\hfill $\diamond$\medskip}

\begin{document}
\sloppy
\maketitle
\begin{abstract} In this paper, we study the structural stability of three-dimensional diffeomorphisms with source-sink dynamics. Here the role of source and sink is played by one-dimensional hyperbolic repeller and attractor. It is well known that in the case when the repeller and the attractor are solenoids (not embedded in the surface), the diffeomorphism is not structurally stable. The author proves that in the case when the attractor and the repeller are canonically embedded in a surface, the diffeomorphism is also not structurally stable.
\end{abstract}
\section{Introduction and formulation of results}
The dynamics of any $\Omega$-stable diffeomorphism $f$ of a closed connected $n$-manifold\footnote{All the manifolds considered in the paper are assumed to be orientable and all diffeomorphisms are assumed to preserve orientation.} $M^n$ can be represented as an attractor-repeller. In this situation, all points outside the attractor and the repeller are wandering and move from the repeller to the attractor\footnote{A set $A$ is called an {\it attractor} of a diffeomorphism $f$ if it has a compact neighborhood $ U_A$ such that $ f (U_A) \subset \textit {int} ~ U_A$ and 
$A= \bigcap \limits_{k \geqslant 0} f^k {(U_A)} $.  $ U_A$ is called a {\it trapping neighborhood} of $A$,  $\bigcup \limits_{k \in\mathbb Z} f^k {(U_A)}$ is called a {\it basin} of the attractor $A$. A {\it repeller} is defined as the  attractor for $ f^{- 1} $. By a {\it dimension} of the attractor (repeller) we mean its topological dimension.}. Such a representation is not unique in the general case and is a source of finding topological invariants of both regular and chaotic dynamical systems. If the attractor and repeller are basic sets, the wandering set is foliated by stable attractor manifolds and unstable repeller manifolds simultaneously. If the topological dimensions of the attractor and the repeller coincide, then the foliations also have the same dimension. The subject of many studies is the question on the existence and the structural stability of such a diffeomorphism.

\begin{figure}[h!]
\centerline{\includegraphics
[width=5cm,height=5cm]{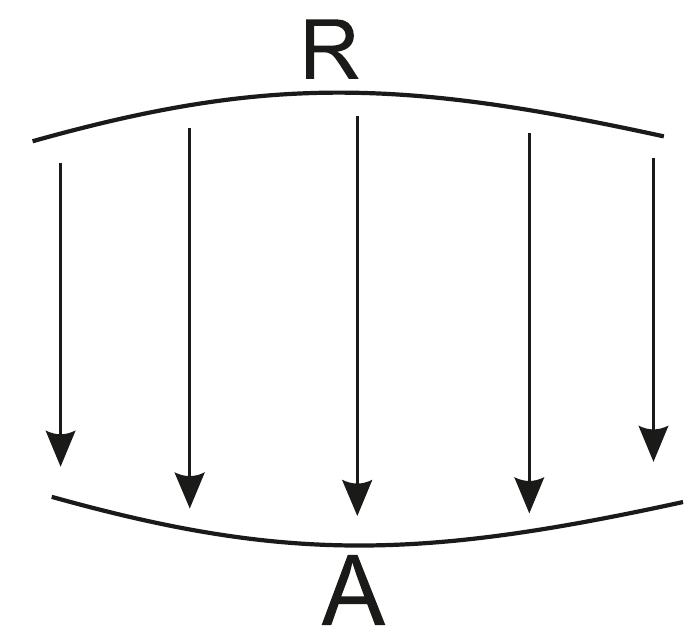}}
\caption{Dynamic attractor-repeller}\label{ARglo}
\end{figure}

Thus, we are dealing with the following situation (see Fig. \ref{ARglo}).
\begin{itemize}
\item $M^n$ is a smooth closed $n$-manifold 
\item $f:M^n\to M^n$ is a preserving orientation $\Omega$-stable diffeomorphism whose non-wandering set consists of two basic sets, an attractor $A$ and a repeller $R$  of the same topological dimension $k\in\{0,\dots,n-1\}$\footnote{Since the topological dimension of the basic set does not exceed the dimension of the supporting manifold, $k$ can take values from $0$ to $n$. However, in the case of $k=n$, the basic set is unique, coincides with the ambient manifold and $f$ is an Anosov diffeomorphism (see, for example, \cite[Theorem 8.1]{GrPo}).}. 
\end{itemize} 

\begin{figure}[h!]
\centerline{\includegraphics
[width=5cm,height=5.6cm]{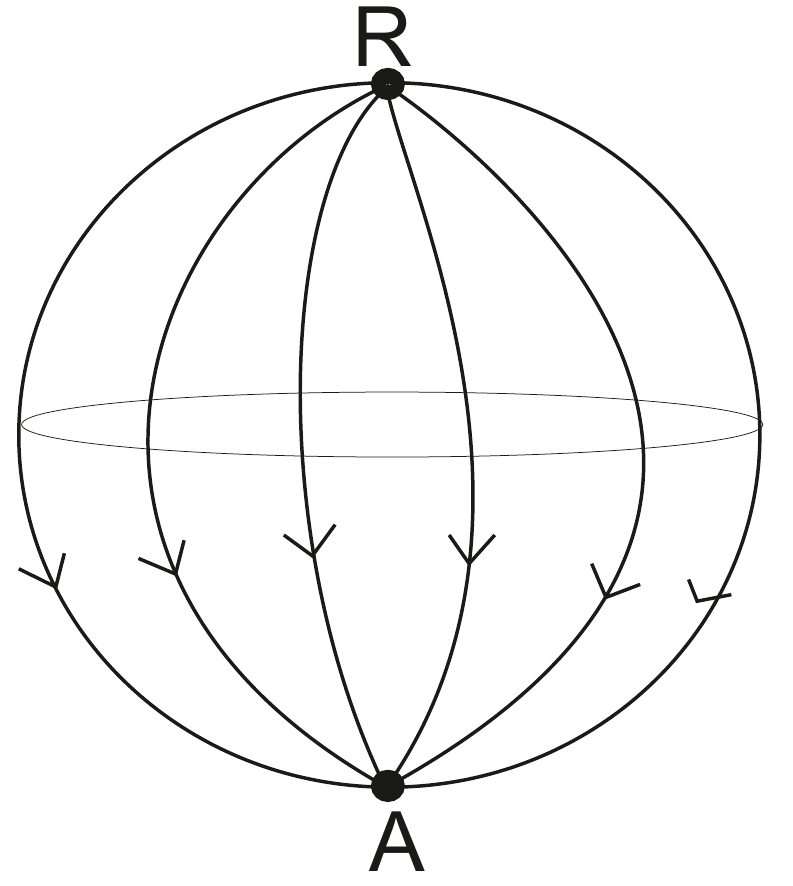}}\caption{Sink-source diffeomorphism}\label{NS-sn}
\end{figure}

If $\bf k=0$, then the attractor and the repeller are sink and source, respectively (see Fig. \ref{NS-sn}). In this case, all such diffeomorphisms are given on the $n$-sphere,  are structurally stable, and form a unique class of topological conjugacy (see, for example, \cite[Theorem 2.5]{GrPo}).

\begin{figure}[h!]
\centerline{\includegraphics
[width=8cm,height=8cm]{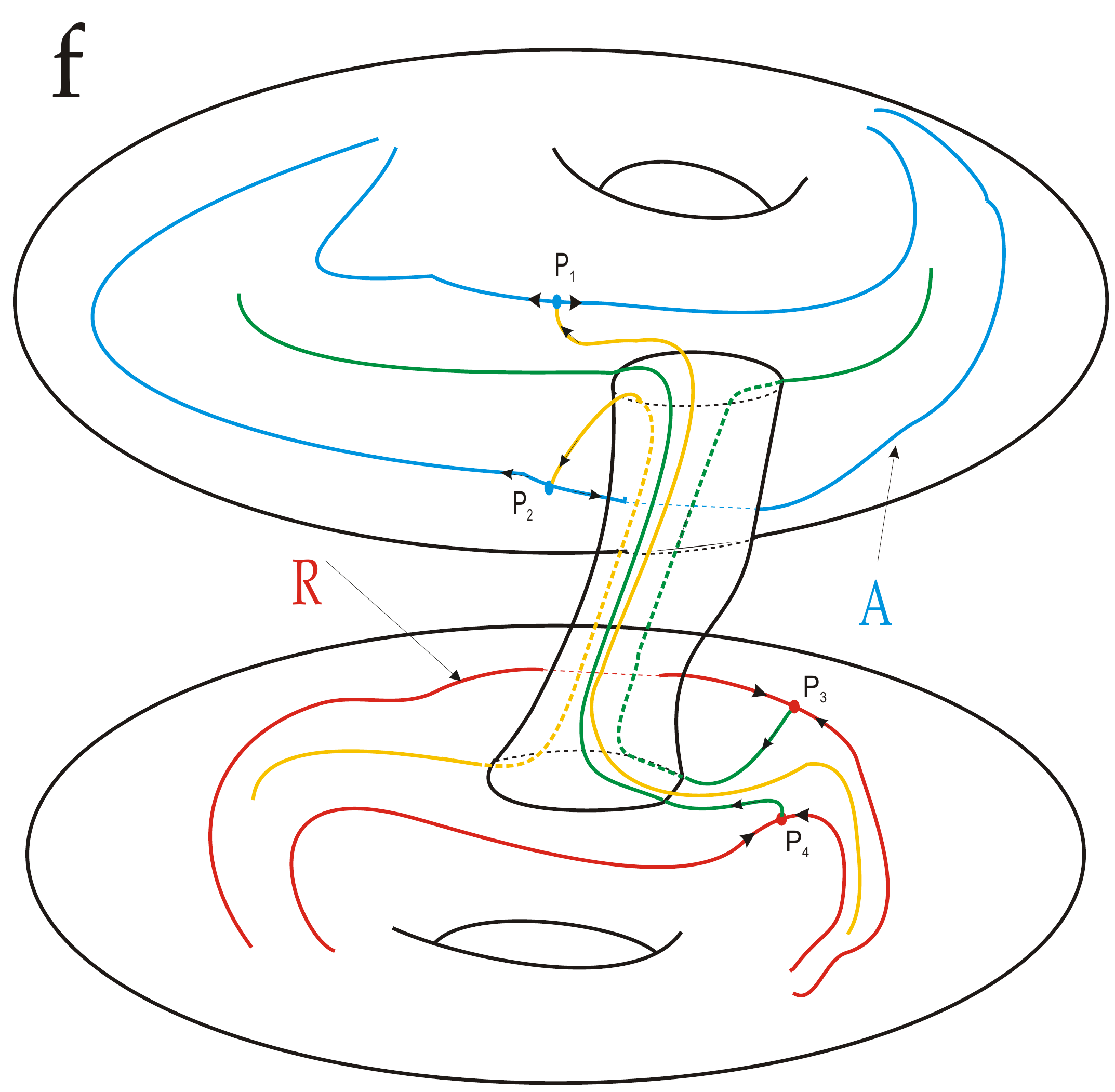}}\caption{Robinson-Williams example}\label{AR-RW}
\end{figure}

If $\bf k=1,\,n=2$ the dynamics on the attractor and repeller are conjugated with the dynamics on non-trivial basic sets of the so-called DA (derived from an Anosov) or DPA (derived  from a pseudo-Anosov) surface diffeomorphisms (see, for example, \cite{PlZh}).
Particular case is a surfaces diffeomorphism with one-dimensional attractor and repeller. Such dynamics is achieved, for example, by taking of the connected sum of two DA-models on 2-tori (see Fig. \ref{AR-RW}). R.C. Robinson and R.F. Williams \cite{RoWi} construct an open set of diffeomorphisms with such phase portraits no one of which is finitely stable.

It follows from \cite{GrMi} that there are no structurally  stable $2$-diffeomorphisms whose non-wandering set is a disjoint union of one-dimensional hyperbolic attractor and repeller.

\begin{figure}[h]\centerline{\includegraphics[width=0.6\textwidth]{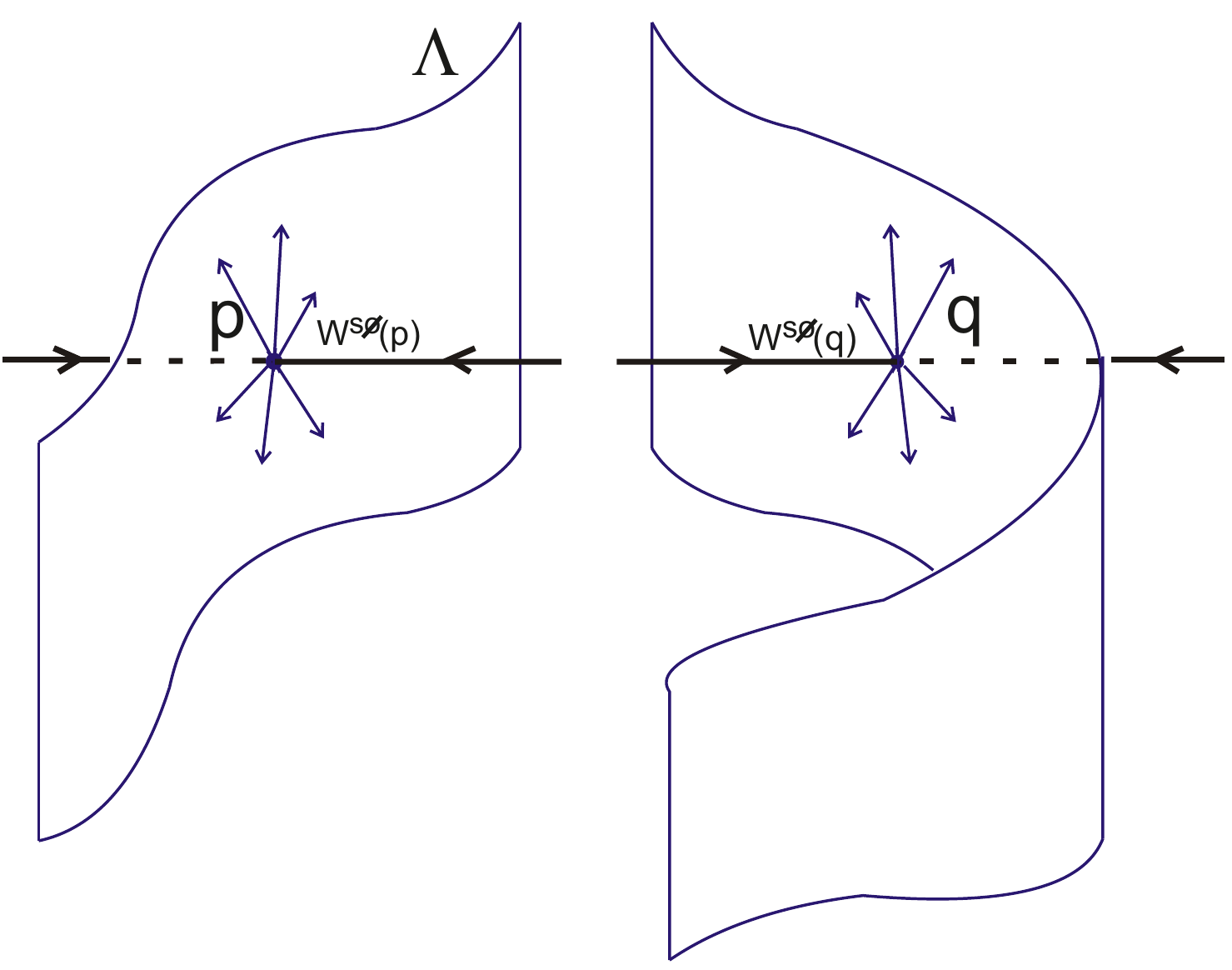}}\caption{Expanding 2-dimensional attractor, obtained from Anosov diffeomorphism by the surgery operation}\label{23-bunch}
\end{figure}

If $\bf k=2,\,n=3$ then due to Brown result \cite{Brown2010} the attractor (the repeller)  is either expanding attractor (contracting repeller) or surface attractor (surface repeller). Recall that, due to \cite{Wi74}, an attractor $A$ of $f$ is said to be {\it expanding} if the topological dimension of   $A$  is equal to  the dimension of $W^u_x, x\in A$ (see Fig. \ref{23-bunch}). One says that $R$ is a  {\it  contracting repeller} of  $f$ if it is  an expanding attractor for $f^{-1}$.
According to \cite{GrMeZh}, a hyperbolic attractor of a diffeomorphism $f:M^3\to M^3$ is called a {\it surface attractor} if it is contained in a compact  surface (not necessarily connected and possible with boundary) $\Sigma_A$ topologically embedded in $M^3$ such that $f(\Sigma_A)\subset \Sigma_A$. One says that $R$ is a  {\it  surface repeller} of  $f$ if it is  a surface attractor for $f^{-1}$.  

\begin{figure}[h!]
\centerline{\includegraphics
[width=12cm,height=6cm]{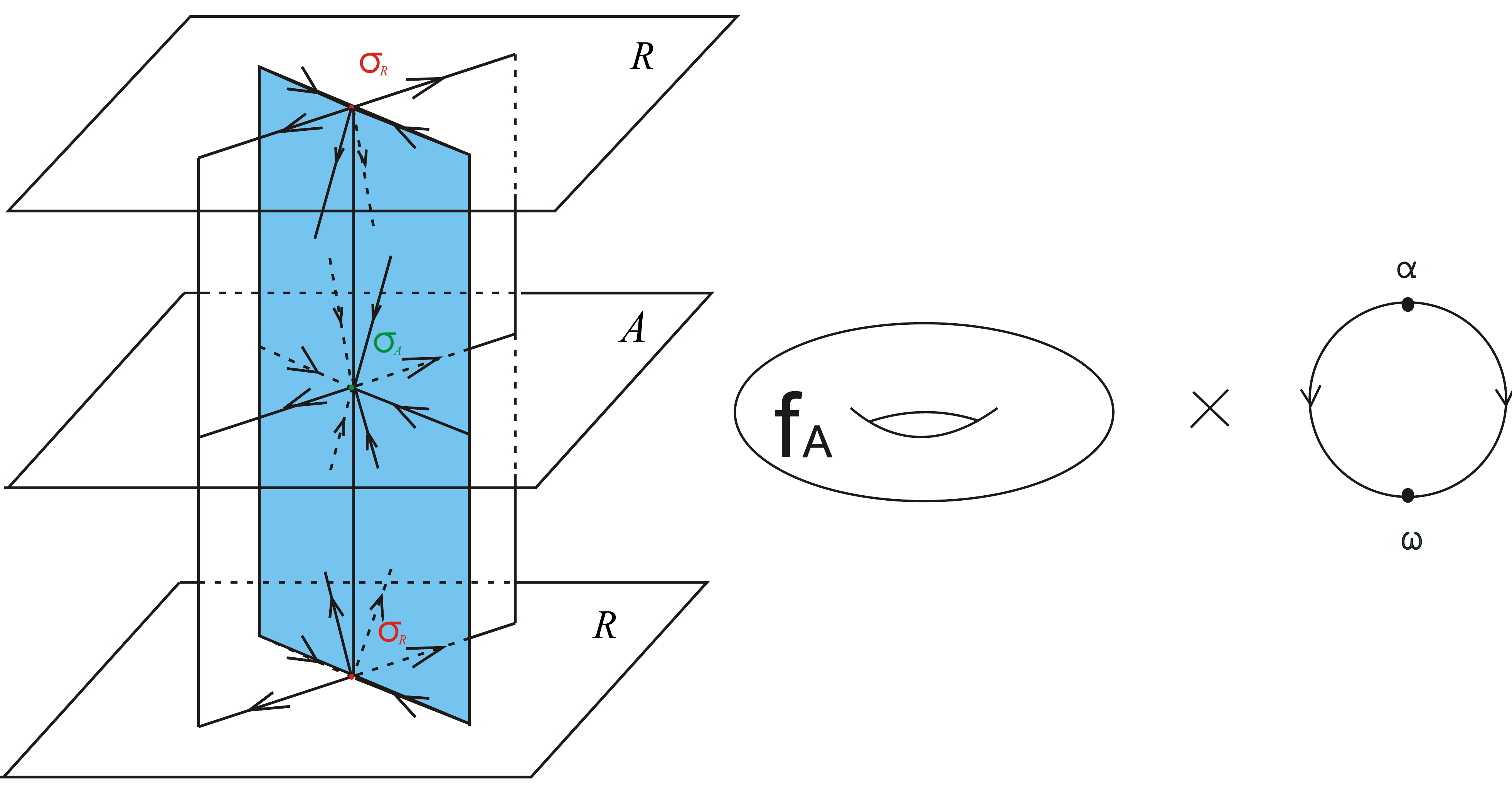}}\caption{Structural stable 3-diffeomorphism with 2-dimensional surfaced attractor and repeller}\label{G}
\end{figure}

It follows from results by V. Grines, V. Medvedev, E. Zhuzhoma  \cite{GrinesZh2005},  \cite{MeZhu2005} that the dynamics is always not structural stable if either $A$ is an expanding attractor or $R$ is a contracting repeller. V. Grines, Yu. Levchenko, V. Medvedev, O. Pochinka \cite{GrLePo} proved that dynamics where both the attractor and the repeller are surface there is only on mapping torus, it possible to be structural stable (see Fig. \ref{G}) and the authors obtained complete topological classification of such rough systems.     

\begin{figure}[h!]
\centerline{\includegraphics
[width=6cm,height=4.5cm]{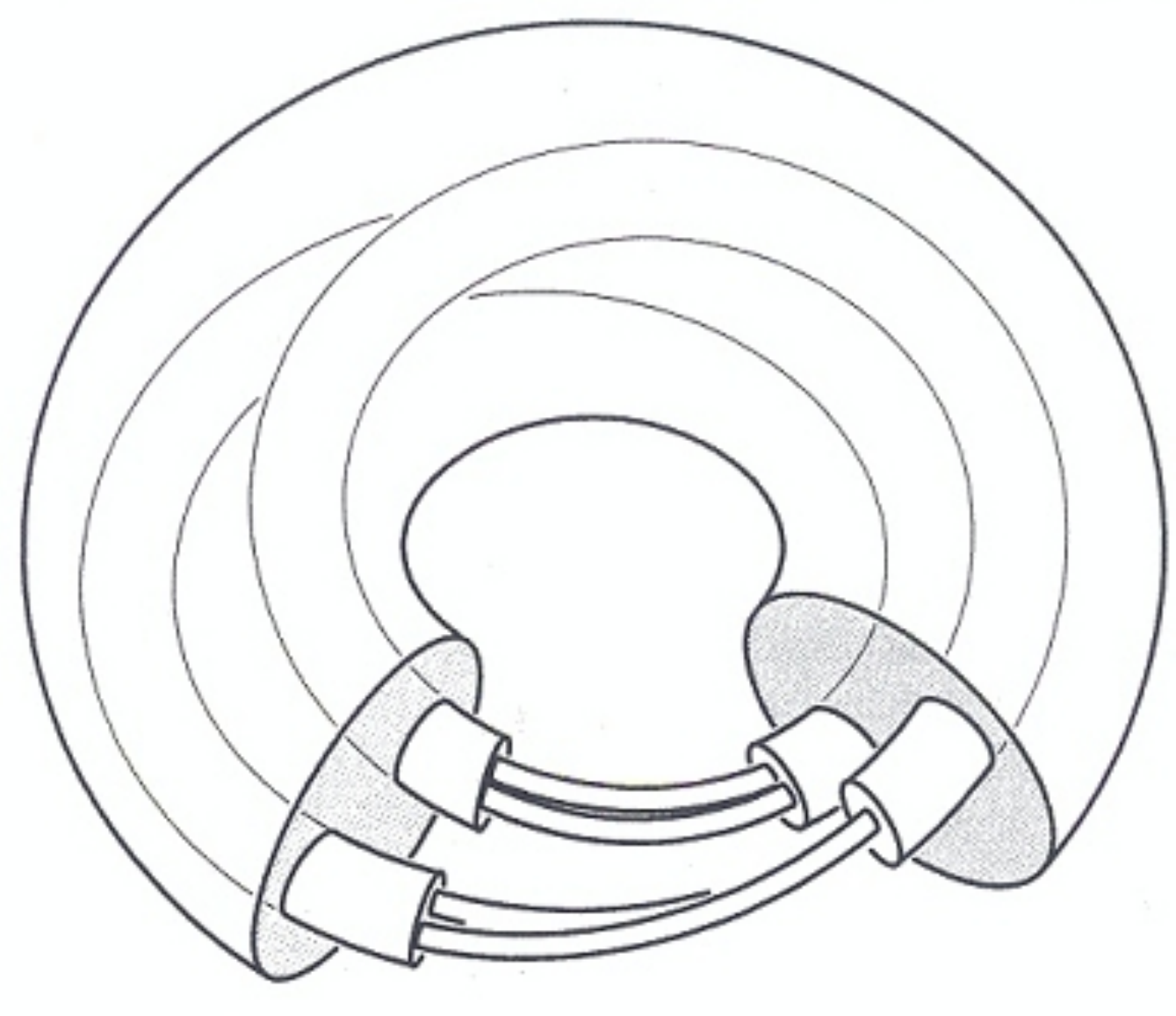}}\caption{Smale's solenoid}\label{solenoid}
\end{figure}

If $\bf k=1,\,n=3$ then the attractor (the repeller) is automatically expanding (contracting) as it consists of the unstable (stable)
manifolds of its points, that was proved by R. Plykin \cite{Plykin1971}.
R. Williams \cite{Wi74} shows that the dynamics on such a basic set  is conjugate to the shift on the 
reverse limit of a branched 1-manifold with respect to an expanding map. A construction of 3-diffeomorphisms with one-dimensional attractor-repeller dynamics firstly was suggested by J. Gibbons  \cite{Gib72}. He construct many models on 3-sphere with Smale's solenoid basic sets (see Fig. \ref{solenoid}) and proves that all examples are not structurally stable. B. Jiang, Y. Ni and S. Wang \cite{JNW} proved that a 3-manifold $M^3$ admits a  diffeomorphism $f$ whose non-wandering set consists of Smale's solenoid attractors and repellers if and only if $M^3$ is a lens space $L(p, q)$ with $p\neq 0$. They also shown that such $f$ are not structural stable. 

All generalizations of  Smale's solenoid as the intersections of nested handlebodies are not surface. Moreover, all known examples of diffeomorphisms with the generalized solenoids as the attractor and the repeller are not structurally stable.

\begin{figure}[h]
\centerline{\includegraphics[width=8 true cm]{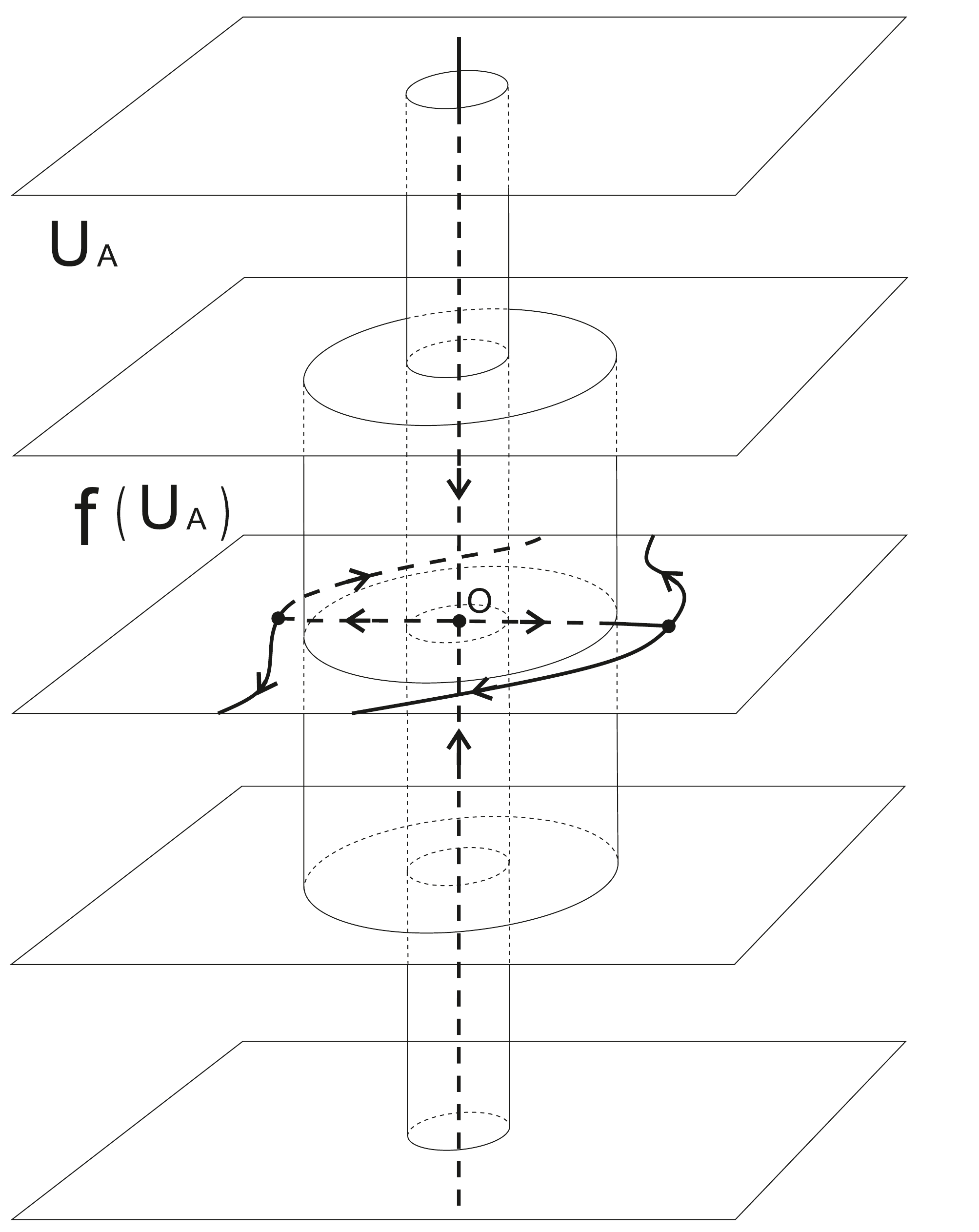}}
\caption{Trapping neighborhood of a canonically embedded surface attractor $A$}\label{ne}
\end{figure}

A natural way to get a surface one-dimensional attractor for a 3-diffeomorphism $f$ is to take an attractor $A$ of some 2-diffeomorphism and multiply its trapping neithborhood by a contraction in transversal direction (see Fig. \ref{ne}). According to \cite{BGPY} such attractor $A$ is called {\it canonically embedded surface attractor}. One says that $R$ is a  {\it canonically embedded surface repeller} of  $f$ if it is  a surface attractor for $f^{-1}$. 

Infinitely many pairwise $\Omega$-non-conjugated diffeomorphisms with such attractors and repellers  were constructed in \cite{BGPY}. Moreover, there a conjecture was formulated that all such  diffeomorphisms are not structurally stable. The main result of this paper is the proof of the conjecture.
\begin{theo}\label{the-hyp} There are no structurally  stable $3$-diffeomorphisms whose non-wandering set is a disjoint union of one-dimensional hyperbolic canonically embedded surface attractor and repeller. 
\end{theo}

Notice, that in \cite{Bonatti2010}, \cite{Shi2014} structurally stable 3-diffeomorphisms with one-dimensional attractor-repeller  dynamics were constructed, but  the constructed basic sets were  not canonically embedded in  surfaces in that examples.

\section{One dimensional basic sets for diffeomorphisms of surfaces}\label{bss}
Let $M^2$ be a closed surface and $\psi : M^2 \to M^2$ be an $\Omega$-stable diffeomorphism. In this section we describe important properties of
one-dimensional basic sets for diffeomorphisms of surfaces following by \cite{Plykin1971},   \cite{Grines1975} (see also \cite{GrPo}). 

For simplicity everywhere below we assume that the basic set is a  connected attractor $A$. Then  
\begin{itemize}
\item  $ A = \bigcup \limits_{x \in A} W^u_x $;
\item at least one of the connected components of the set $ W^s_x \setminus \{x \}, x \in A $ contains a dense set in $ A $;
\item there are a finite number of points $ x \in A $ for which one of the connected components $W^{s-}_x$ of the set $ W^s_x \setminus \{x \} $ does not intersect $ A $, denote by $W^{s+}_x$ other connected component. Such points are called {\it $ s $-boundary}, their set $P_A$ is not empty and consists of a finite number of periodic points;
\item the set $W^s_{A}\setminus A$ consists of a finite number path-connected components. 
\end{itemize}

\begin{figure}[h]
\centerline{\includegraphics[width=8 true cm]{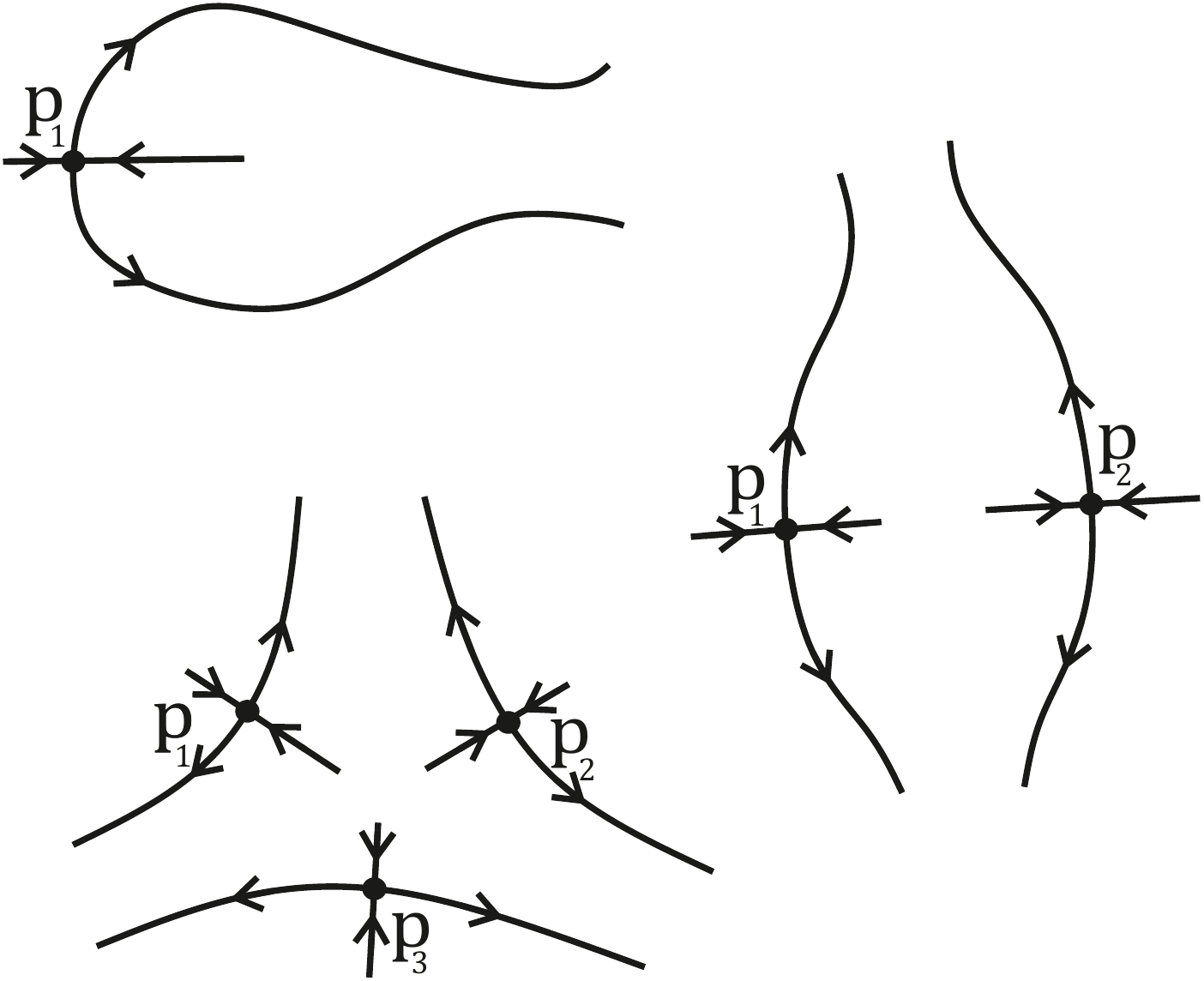}}
\caption{Bunches of degrees 1, 2 and 3}\label{periodical_points}
\end{figure}

A {\it bunch $b$} of the attractor $A$ is the union of the maximal number $r_b$ of the unstable manifolds $W^{u}_{p_{1}},\dots,W^{u}_{p_{r_b}}$ of the $s$-boundary points $p_{1}, \dots,
p_{r_b}$ of the set $A$ for which $W^{s-}_{p_1},\dots,W^{s-}_{p_{r_b}}$ belong to the same path-connected components of $W^s_{A}\setminus A$. The number $ r_b $ is called a {\it degree of the bunch} (see Fig. \ref{periodical_points}). Let $ B _{A} $ be the set of all bunches of the attractor $ {A} $. 

\begin{figure}[h]
\centerline{\includegraphics[width
=8 true cm]{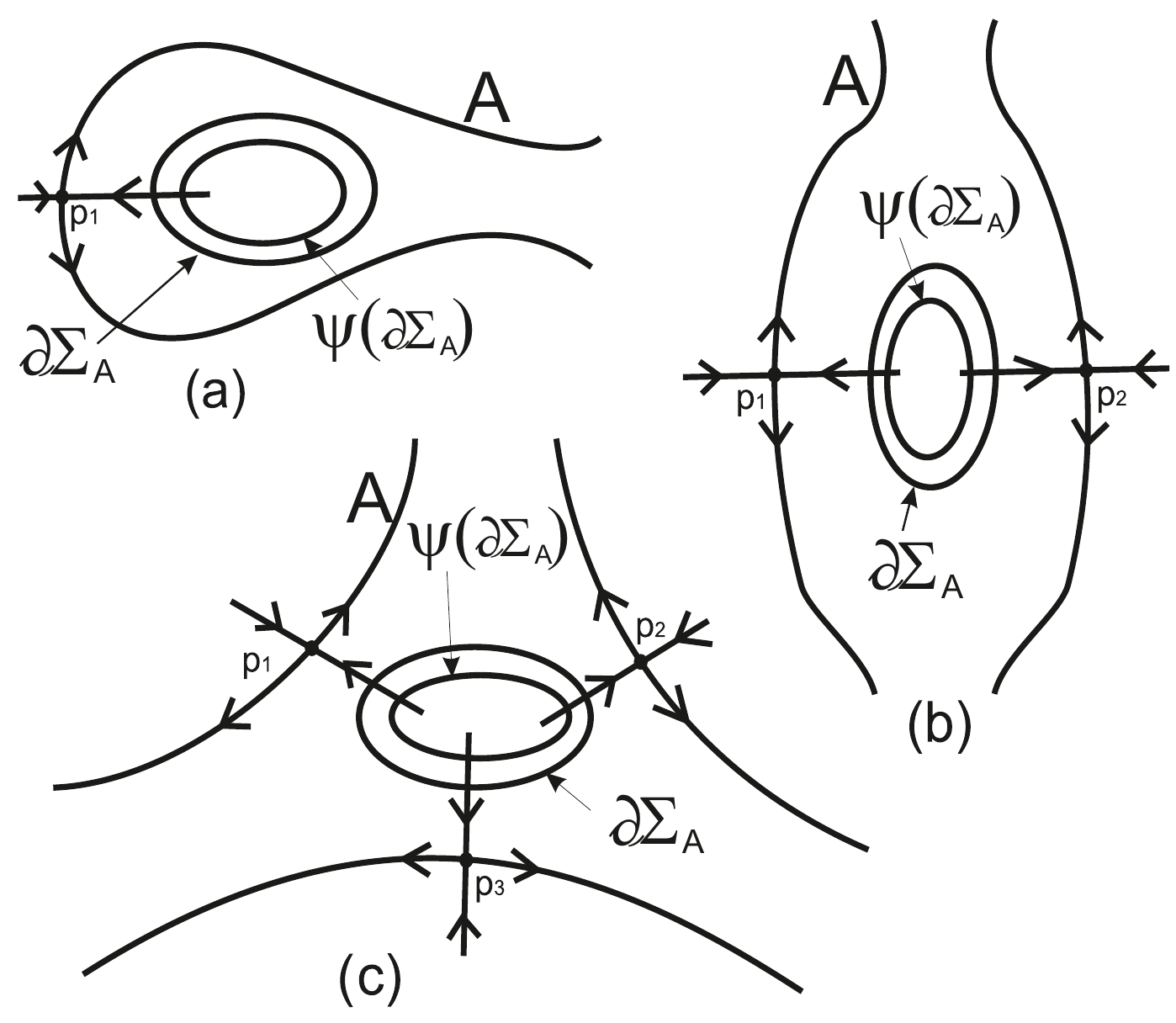}}
\caption{Canonical supports of 1-dimensional attractors}\label{11}
\end{figure}

For attractor $A$ let $m_{A}$ denote the number of its bunches and let $r_{A}$ denote the sum of the degrees of these bunches.  Then  
\begin{itemize}   
\item  there is a trapping neighborhood $\Sigma_A$ ($\psi(\Sigma_A)\subset int\,\Sigma_A$) of the attractor $A$ such that $\Sigma_{A}$ is a compact orientable surface of genus $g_{A}=1+\frac{r_{A}}{4}-\frac{m_{A}}{2}$, it has $m_{A}$ boundary components and it has negative Euler characteristic (see Fig. \ref{11}).
\end{itemize} 

\section{Point of view from Lobachevsky plane}\label{Lob}
In this section, we describe the dynamics on a one-dimensional surface basic set from the point of view of its lifting to the Lobachevsky plane.
\subsection{Properties of automorphisms of deck transformation groups}\label{desc}
In this section we describe important properties of automorphisms of deck transformation groups for universal cover of surfaces following by \cite{HTh} (see also \cite{GrPo}).

We consider the Poincar\'{e} disk model of the hyperbolic plane  as the unit open ball $\mathbb U=\{z\in\mathbb C~:~\vert z
\vert < 1\}$ of the complex plane with the hyperbolic metric $d$. The boundary of the ball $\mathbb U$  is called the {\it absolute of the hyperbolic plane} denoted by $\mathbb E$ ($\mathbb E=\partial{\mathbb
U}=\{z\in\mathbb C~:~\vert z
\vert = 1\}$). 

Let us describe a characteristic of so called {\it hyperbolic isometry} $g:\mathbb U\to\mathbb U$. 
\begin{itemize}
\item $g$ uniquely extends to the absolute $\mathbb E$ and has exactly two fixed points in $\mathbb E$ and it has no fixed points in $\mathbb U$; 
\item there is the unique geodesic $l_g$ which is invariant under $g$. This geodesic is called the {\em axis of the hyperbolic isometry $g$}. The axis joins the points $P_g$ and $Q_g$ on the absolute. The restriction of $g$ to its axis is a shift, i.e. $d(x,g(x))=d_g$ for every point $x\in l_g$ (see Fig. \ref{lg});
\begin{figure}[h]
\centerline{\includegraphics [width = 7 true cm]{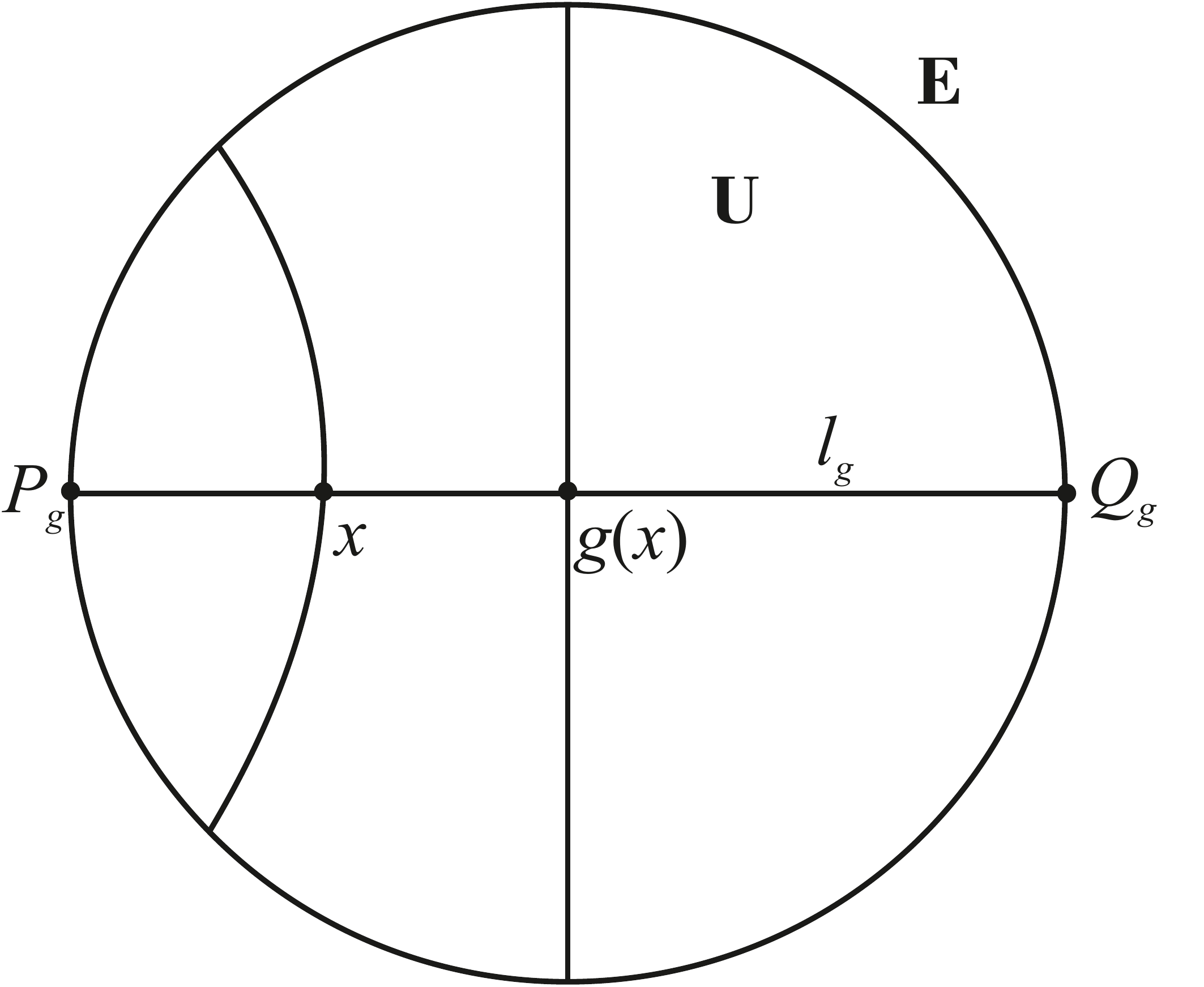}}
\caption{Action of a hyperbolic isometry on the Poincar\'{e} disk} \label{lg}
\end{figure}
\item two distinct hyperbolic isometries $g_1,g_2$  have a common fixed point if and only if there is a hyperbolic isometry  $\gamma$ and there are integers $k_1,k_2$ such that $g_1=\gamma^{k_1}$ and $g_2=\gamma^{k_2}$. Therefore, two hyperbolic isometries $g_1,g_2$ either have no common fixed points or both fixed points of $g_1$ are the fixed points of $g_2$ as well. 
\end{itemize}
We say the fixed points of each hyperbolic isometry to be {\it rational} and call {\it irrational} any other point of the absolute. 

\begin{figure}[h]
\centerline{\includegraphics [width = 12 true cm]{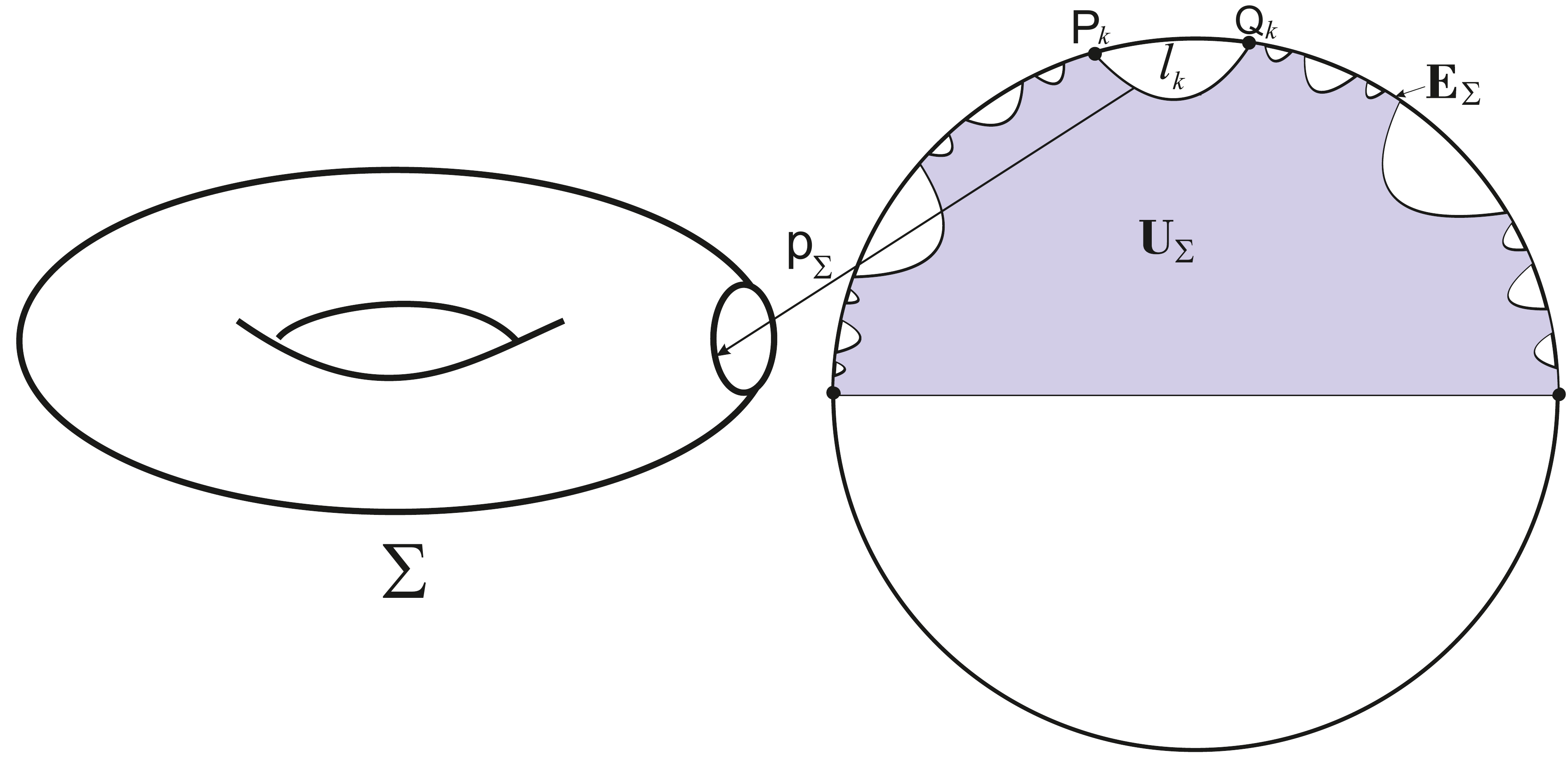}}
\caption{A universal cover for  $\Sigma$} \label{eba}
\end{figure}

Now let $\Sigma$ be an orientable surface with boundary (possible empty) of negative Euler characteristic. If the boundary is empty then we put $S=\Sigma$, in the opposite case we glue two copies of $\Sigma$ along the boundary components we get a surface $S$ without boundary. The curves along which we glue are essential and, therefore, they can be assumed to be geodesics. Then
\begin{itemize}
\item there is a subgroup $\mathbb G_{\Sigma}$ of the hyperbolic isometrics of $\mathbb U$,  which is isomorphic to the fundamental group of the manifold $\Sigma$ and  a connected set $\mathbb U_{\Sigma}\subset\mathbb U$, which universally covers $\Sigma$. That is $\Sigma=\mathbb U_\Sigma/\mathbb G_{\Sigma}$ and $p_{\Sigma}:\mathbb U_\Sigma\to \Sigma$ is a universal cover. 
\item the set $\mathbb Q_\Sigma$ of the rational points of the isometries from $\mathbb G_{\Sigma}$ is countable and dense on ${\mathbb E}_{\Sigma}=\partial\mathbb U_\Sigma\cap\mathbb E$;
\item $\mathbb E_{\Sigma}$ is the Cantor perfect set on the absolute and, therefore, it can be expressed as  $\mathbb
E\setminus\bigcup\limits_{k \in {\mathbb N}}(P_{k},Q_{k})$ where $(P_{k}, Q_{k})$ are the adjacent intervals of the Cantor set $\mathbb E_{\Sigma}$. Every geodesic $l_{k}\subset {\mathbb U}$ with the boundary points $P_{k}, Q_{k}$ belongs to $\mathbb U_{\Sigma}$, is the axis of some non-identity element $g\in\mathbb G_{\Sigma}$ and $p_\Sigma(l_k)$ is a connected component of $\partial\Sigma$  (see Fig. \ref{eba}). 
\end{itemize}
  
Every automorphism $\tau: \mathbb G_{\Sigma}\to \mathbb G_{\Sigma}$ induces on the set $\mathbb Q_\Sigma$ the map which sends the fixed points of an element $g$ onto the fixed points of the element $\tau(g)$. This map uniquely extends to the homeomorphism $\tau^e:\mathbb E_\Sigma\to\mathbb E_\Sigma$. 
Also every automorphism $\gamma: \mathbb G_{\Sigma}\to \mathbb G_{\Sigma}$ induces the {\it internal automorphism} $A_\gamma: \mathbb G_{\Sigma}\to \mathbb G_{\Sigma}$ by the formula $$A_\gamma(g)=\gamma g\gamma^{-1},\,g\in\mathbb G_{\Sigma}.$$
Let $\mathcal A_\tau=\{A_{\gamma}\tau^k,~\gamma\in\mathbb G_{\Sigma},~k\in\mathbb Z\}$.
An automorphism $\tau$ of the group $\mathbb G_{\Sigma}$ is said to be {\it hyperbolic} if  $a(g)\neq g$ for every non-identity automorphism $a\in\mathcal A_\tau$ and for every non-identity element $g$ except those that  $p_{\Sigma}(l_{g})\subset\partial\Sigma$. 

Let  $h:\Sigma\to\Sigma$ be a homeomorphism. Since $\mathbb U_\Sigma$ is the universal covering space, it implies that there is a homeomorphism (not unique) $\bar h: \mathbb U_\Sigma\to\mathbb U_\Sigma$ which is a lift of a map $h$ (i.e. the homeomorphism $\bar h$ satisfies $p_{\Sigma}\bar h=h p_{\Sigma}$). Let $h(x)=y$ and $\bar x\in p^{-1}_{\Sigma}(x)$, $\bar y\in p^{-1}_{\Sigma}(y)$, $\bar h(\bar x)=\bar y$ and $g\in\mathbb G_{\Sigma}$. Since $p_{\Sigma}(g(\bar x))=x$ we have $p_{\Sigma}(\bar h(g(\bar x)))=y$. Hence, $\bar h(g(\bar x))=g'(\bar y)$ for some element $g'\in\mathbb G_{\Sigma}$ and, therefore, $\bar h(g(\bar x))=g'(\bar h(\bar x))$. Thus, the lift $\bar h$ induces the automorphism  $\tau_{\bar h}$ of the group $\mathbb G_{\Sigma}$ which assigns the element $g'$ to an element $g$ by $\tau_{\bar h}(g)=\bar hg\bar h^{-1}$. 
Every lift $\bar h: \mathbb U_\Sigma\to\mathbb U_\Sigma$ of a homeomorphism $h:\Sigma\to\Sigma$ possesses to following properties:
\begin{itemize} 
\item  $\bar h$ uniquely extends onto $\mathbb E_\Sigma$ by the homeomorphism $\bar h^*:\mathbb E\to\mathbb E$ and  $\bar h^*=\tau^e_{\bar h}$;
\item if $h_1,~h:\Sigma\to\Sigma$ are homotopic homeomorphisms then there is a lift $\bar h_1: \mathbb U_\Sigma\to\mathbb U_\Sigma$ of $h_1$ such that $\bar h_1^*=\bar h^*$. 
\end{itemize}
Homeomorphisms $h_1,\,h_2:\Sigma\to\Sigma$ are called {\it $\pi_1$-conjugate} if there are their lifts $\bar h_1,\,\bar h_2: \mathbb U_\Sigma\to\mathbb U_\Sigma$ and an automorphism $\tau:\mathbb G_{\Sigma}\to\mathbb G_{\Sigma}$ such that $\tau\tau_{\bar h_1}=\tau_{\bar h_2}\tau$.

\subsection{Lifting of one-dimensional attractor to the Lobachevsky plane}\label{ALob}
In this section we describe important properties of
one-dimensional basic sets for diffeomorphisms of surfaces following by \cite{Plykin1971},   \cite{Grines1975} (see also \cite{GrPo}). 

Let $A$ be a hyperbolic one-dimensional attractor of an $\Omega$-stable diffeomorphism $f:M^2\to M^2$ and $\Sigma_A$ be its trapping neighborhood (see section \ref{bss}).
Let $p_{\Sigma_{A}}:\mathbb U_{\Sigma_{A}}\to \Sigma_{A}$ be a universal covering and let $\mathbb G_{\Sigma_{A}}$ be the group of its covering transformations. According to the section \ref{desc} every lifting $\bar \psi_{_A}$ induces an automorphism $\tau_{{\bar \psi_{_A}}}$ of the group $\mathbb G_{\Sigma_A}$. Then 
\begin{itemize}
\item the automorphism $\tau_{\bar\psi_{_A}}$ is hyperbolic. 
\end{itemize}

Let $\bar A=p_{\Sigma_{A}}^{-1}(A)$. If $a\in A$ then let $\bar a\in\bar A $ denote the point in the preimage $p_{\Sigma_{A}}^{-1}(a)$. Let $\delta\in\{u,s\}$ and $\nu\in\{+,-\}$. Denote by $w_{\bar a}^{\delta}$ the curve on $\mathbb U_{\Sigma_{A}}$ such that
$p_{\Sigma_{A}}(w_{\bar a}^{\delta})=W^{\delta}_a$.  If $t\in \mathbb R$ is a parameter on the curve $W_{a}^{\delta}$  such that $W_{a}^{\delta}(0)=a$ then
$w_{\bar a}^{\delta}(t)$ is the point on $w_{\bar a}^{\delta}$ such that $p_{\Sigma_{A}}(w_{\bar
a}^{\delta}(t))=W_{a}^{\delta}(t)$ and $w_{\bar a}^{{\delta}+}$, $w_{\bar a}^{{\delta}-}$ are the connected components of the curve $w_{\bar a}^{{\delta}}\setminus \bar a$ for $t>0$, $t<0$ respectively. 

Let $a\in A$. We say that a curve $w_{\bar a}^{{\delta}{\nu}}$ {\it has the asymptotic direction ${\delta}_{\bar a}^{{\nu}}$} for $t\to {\nu}\infty$ if the set
$cl\,(w_{\bar a}^{{\delta}{\nu}})\setminus w_{\bar a}^{{\delta}{\nu}}$ consists of the point $\bar a$ and the point ${\delta}_{\bar a}^{{\nu}}$ which belongs to $\mathbb E_{\Sigma_{A}}$. 
For the attractor $A$ the following asymptotic properties take place (see Figure \ref{soversh}):
\begin{figure}[h]
\centering{\includegraphics[width
=14cm,height=7cm]
{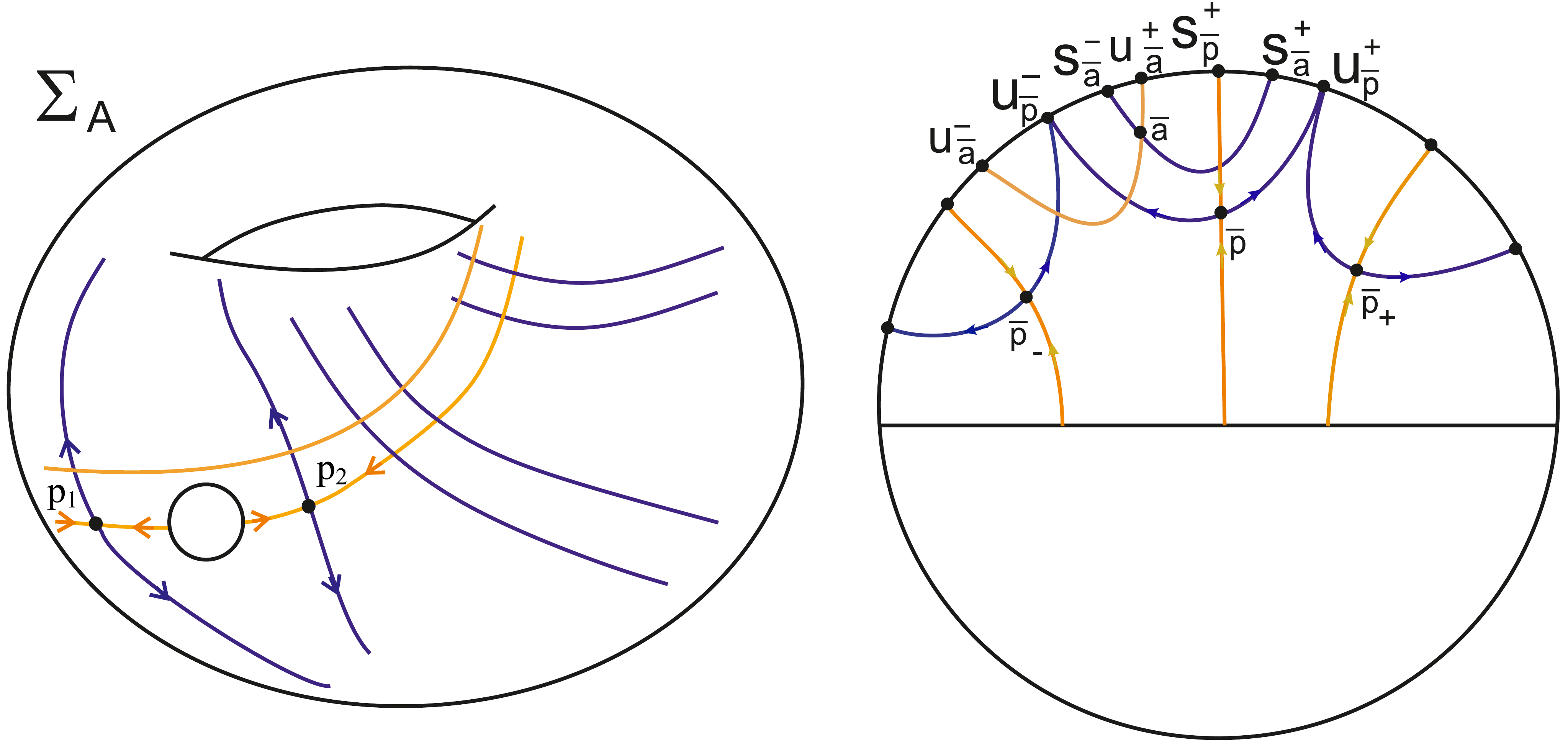}} \caption{Asymptotic behaviour of invariant manifolds on the universal cover, here $p=p_1,\,p_+=p_-=p_2$}\label{soversh}
\end{figure}  
\begin{itemize}
\item if for a point $a\in A$ the component $W^{{\delta}{\nu}}_a$ is densely situated in $A$ then $w^{{\delta}{\nu}}_{\bar a}$ has an irrational asymptotic direction ${\delta}_{\bar a}^{{\nu}}$ for
$t\to {\nu}\infty$;
\item the curve $w_{\bar a}^{u}$ has two distinct boundary points (asymptotic directions) $u_{\bar a}^{+},~u_{\bar a}^{-}$;
\item the curves $w_{\bar a}^{s}$ for some $\bar a\in(\bar A\setminus w^{s+}_{\bar P_A})$ has two distinct boundary points  $s_{\bar a}^{+},~s_{\bar a}^{-}$;
\item for every point $\bar a\in\bar A$ the intersection $cl\, (w^{u}_{\bar{a}})\cap cl\, (w^{s}_{\bar{a}})$ consists of a unique point $\bar a$; 
\item if $w^{s}_{\bar{a}} \cap
w^{s}_{\bar{b}}=\emptyset$ for some points $\bar a,\bar b\in \bar A$ then $cl\, (w^{s}_{\bar{a}}) \cap cl\,(w^{s}_{\bar{b}})=\emptyset$;
\item if $w^{u}_{\bar{a}} \cap
w^{u}_{\bar{b}}=\emptyset$ for some points $\bar a\in\bar A,~\bar b\in(\bar A\setminus w^u_{\bar P_A})$ then $cl\,(w^{u}_{\bar{a}}) \cap cl\,(w^{u}_{\bar{b}})=\emptyset$;
\item for every point $\bar p\in\bar P_A$  there are two distinct points $\bar p_+,\bar p_-\in \bar P_A$ such that $w^u_{\bar p_+}\cap w^u_{\bar p}=\emptyset$, $w^u_{\bar p_-}\cap w^u_{\bar p}=\emptyset$,  $cl\,(w^u_{\bar p_+})\cap cl\,(w^u_{\bar p})=u^{+}_{\bar p}$, $cl\,(w^u_{\bar p_-})\cap cl\,(w^u_{\bar p})=u^{-}_{\bar p}$ and the points $p,p_-,p_+$ are the $s$-boundary points of the same bunch of the attractor $A$;
\item if  $w^{u}_{\bar{a}} \cap
w^{s}_{\bar{b}}=\emptyset$ for some points $\bar a,\bar b\in\bar A$ then $cl\,(w^{u}_{\bar{a}}) \cap cl\,(w^{s}_{\bar{b}})=\emptyset$.
\end{itemize}

\begin{figure}[h]
\centering{\includegraphics[width
=8.5cm,height=9cm]
{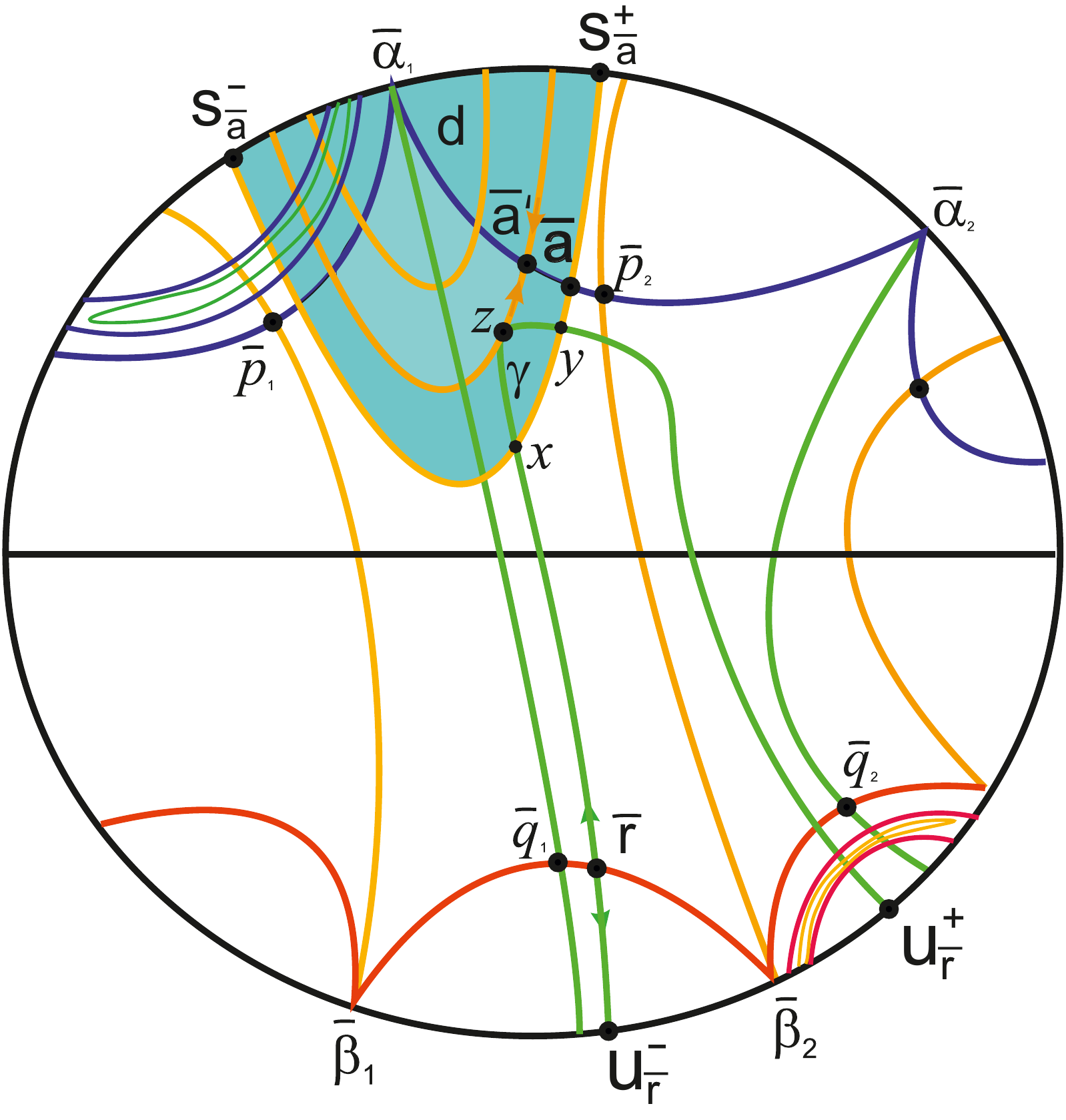}} \caption{The lifting of the Robinson-Williams example on the Lobachevsky plane}\label{RWL}
\end{figure}

At the end of this section we illustrate why the Robinson-Williams example is not structural stable. A phase portrait of the lifting of the Robinson-Williams example on the Lobachevsky plane is given on Fig. \ref{RWL}. It follows from the description above that there are point $\bar a\in\bar A,\,\bar r\in\bar R$ such that $w^{s}_{\bar a}\cap w^{u}_{\bar r}\neq\emptyset$ and $w^{s}_{\bar a},\,w^{u}_{\bar r}$ have asymptotic directions ${s}_{\bar a}^{-},\,{s}_{\bar a}^{+},\,{u}_{\bar r}^{-},\,{u}_{\bar r}^{+}$ which bound disjoint arcs on the absolute. Moreover, the closure of the leaf $w^{s}_{\bar a}$  divides $\mathbb U$ into two connected components, exactly one of them is a 2-disk (denote it $d$) which does not contain the points ${u}_{\bar r}^{-},\,{u}_{\bar r}^{+}$.

If we assume that the leaves $w^{s}_{\bar a}$ and $w^{u}_{\bar r}$ are transversally intersected,  then there is a connected component $\gamma$ of the intersection $w^{u}_{\bar r}\cap d$ which has exactly two intersection points with $w^{s}_{\bar a}$, say $ x, y$. By  Rolle's theorem there is a point  $ z\in\gamma$ where $w^{u}_{\bar r}$ has a contact with a leaf $w^{s}_{\bar a'}$ for a point $\bar a'\in\bar A$, that contradict to the assumption.

\section{Dynamics ``one-dimensional canonically embedded surface attractor-repeller'' for 3-diffeomorphism}\label{dynA-R}

Let $f:M^3\to M^3$ be an $\Omega$-stable whose non-wandering set consists of a one-dimensional canonically embedded surface attractor $A$ and one-dimensional canonically embedded surface repeller $R$. Then  $V_f=M^3\setminus(A\cup R)$ is the wandering set of $f$ which coincides with $W^s_A\setminus A$ and $W^u_R\setminus R$ simultaneously. As our finite aim to refute the structural stability of $f$, everywhere below  we will suppose that all boundary points of the attractor $A$ are fixed (in the opposite case we can consider an appropriate degree of $f$). 

\subsection{The topology of the wandering set}\label{TW}
 By definition of canonically embedded attractor and due to results in section \ref{bss}, $A$ has a trapping  neighborhood $U_{A}$ of the form $\Sigma_A \times [-1,1] $, where $\Sigma_A =\Sigma_A \times \{0 \} $ is a trapping neighborhood of attractor $A$ as an attractor of a 2-diffeomorphism $\psi_{_A}$ and diffeomorphism $ f_{_A}=f|_{U_A}:U_A\to f(U_A)$ has a form  $f_{_A} (w, z) = (\psi_{_A} (w), z/2): \Sigma_A \times [-1,1 ]\to f(\Sigma_A)\times [-1/2,1/2] $ (see Fig. \ref{ne}).  

Recall that $\Sigma_A$ is a compact orientable surface of genus $g_A$ with $m_A$ boundary components and negative Euler characteristic. Let us glue two copies of $\Sigma_A$ along the boundary components we get a closed surface $S_{A}$ of a positive genus $\rho_A=2g_A+m_A-1>1$. Denote by $C_A$ the coinciding copies of $\partial\Sigma_A$ in $S_A$.

\begin{figure}[h]
\centerline{\includegraphics[width
=12 true cm]{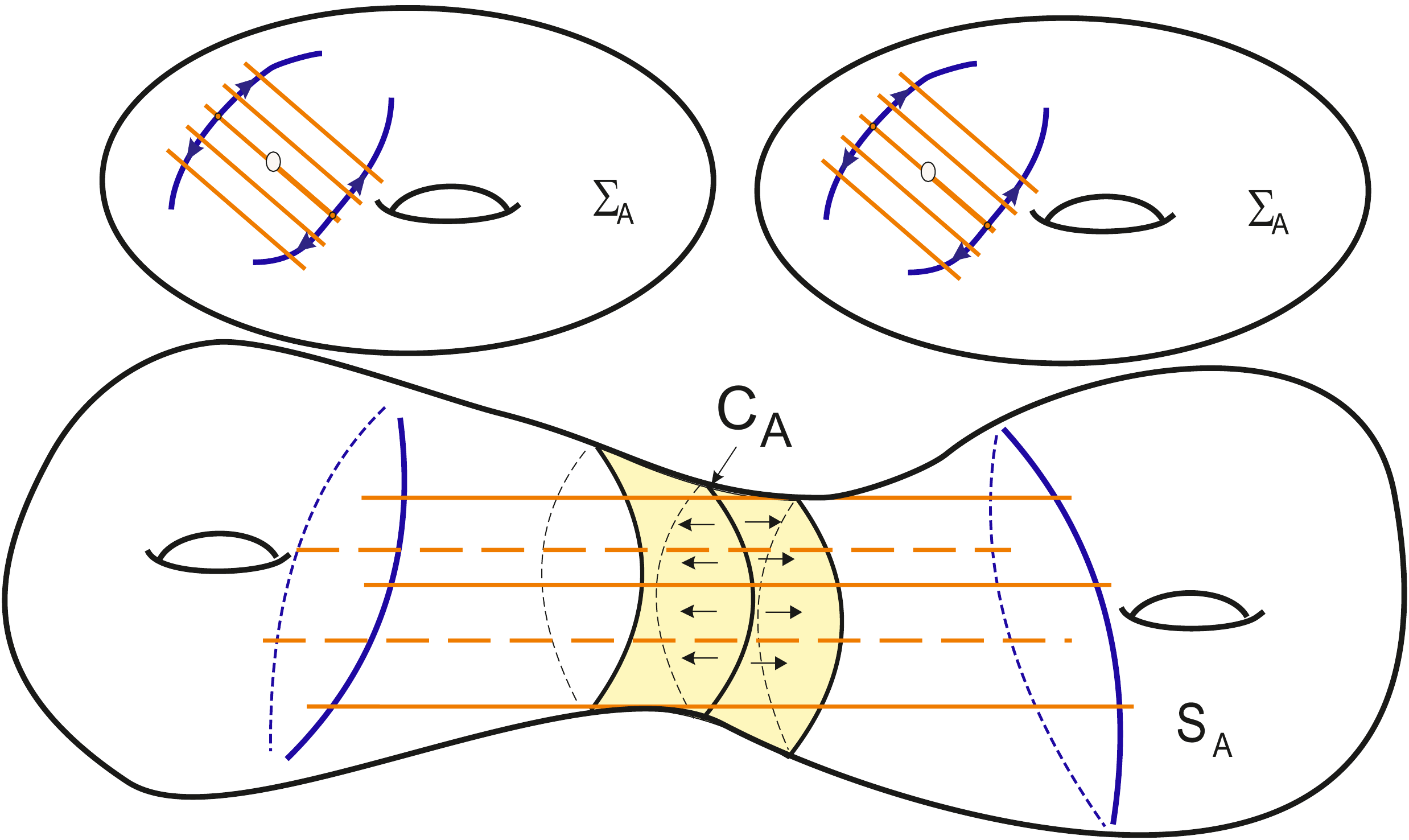}}
\caption{Phase portrait of the diffeomorphism $\Psi_A$}\label{xia}
\end{figure}
\begin{lemm}\label{sia}  Diffeomorphism $f|_{V_f}$ is smoothly conjugate with a diffeomorphism $\chi_A:S_{A}\times\mathbb R\to S_{A}\times\mathbb R$ given by the formula 
$$\chi_A(s,r)=(\Psi_A(s),r+1),$$ where $\Psi_A:S_A\to S_A$ is a reducible diffeomorphism which  coincides with $\psi_{_A}$ on each copy of $\Sigma_A$ outside a neighborhood of $C_A$ and is an extension from $C_A$ in the neighborhood being identical on $C_A$ (see Fig. \ref{xia}).
\end{lemm}
\begin{demo} Let us compactify the surface $\Sigma_A$ by discs $D_A$ to get a closed surface $\tilde\Sigma_A$ and a diffeomorphism $\tilde\psi_{_A}:\tilde\Sigma_{A}\to \tilde\Sigma_{A}$ which is $\psi_{_A}$ on $\Sigma_A$ and every disc in  $D_A$ is a part of the basin of  a linear source point in it. Let $\tilde f_A:\tilde\Sigma_{A}\times\mathbb R\to \tilde\Sigma_{A}\times\mathbb R $ be a diffeomorphism given by the formula $$\tilde f_A (w, z) = (\tilde\psi_{_A} (w), z/2).$$
Let $(x,y)$ be local coordinates in a neighborhood the discs  $D_A$ such that for  $r=\sqrt{x^2+y^2}$ and $\lambda>1$, $D_A=\{(x,y): r\leqslant \lambda^{-4}\}$, $\partial\Sigma_A=\{(x,y): r=1\}$ and 
the diffeomorphism $\tilde\psi_{_A}$  has a form
$$\tilde\psi_{_A}(x,y)=(\lambda x,\lambda  y)$$ for $\lambda^{-2}\leqslant r\leqslant 1$. Let $$L_0=\left\{(r,z): r=\lambda^{-1}-(\lambda^{-1}-\lambda^{-2})\sqrt{1-z^2},\lambda^{-2}\leqslant r\leqslant \lambda^{-1}\right\}.$$  Let $G_0$ coincides with $L_0$ for 
$\lambda^{-2}\leqslant r\leqslant \lambda^{-1}$ and coincides with $M_0=\tilde\Sigma_{A}\times\{-1,1\}$ outside the set $\{(r,z):r<\lambda^{-1}\}$ (see Fig. \ref{signu}). By the construction $G_0\cong S_{A}$ and we will identify their.  
\begin{figure}[h]\centerline{\includegraphics [width = 9 true cm]{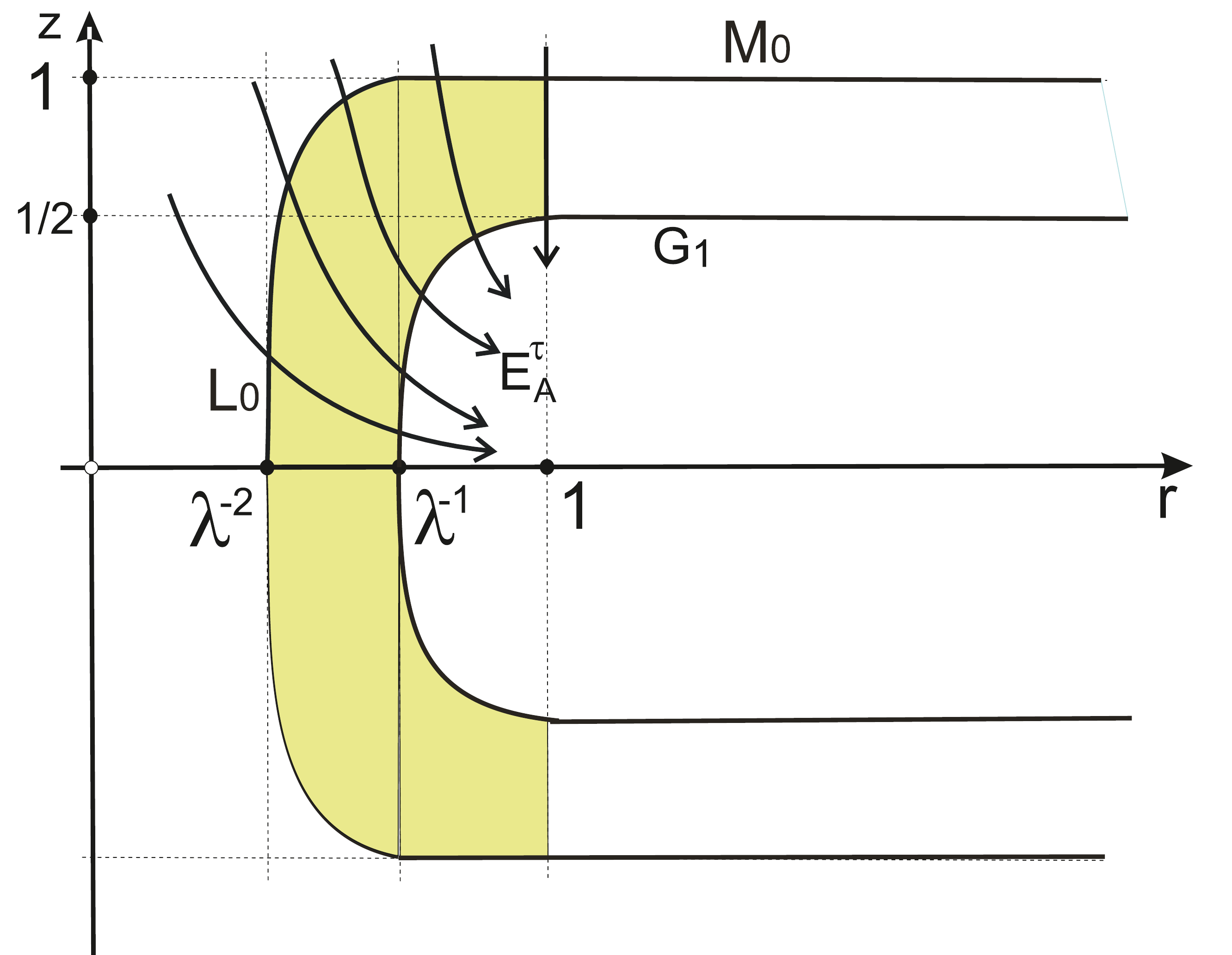}}\caption{Construction of the domain $K_A$}\label{signu}\end{figure}

By the construction the diffeomorphism  $\tilde f_A$ on the set $\{(r,z):\lambda^{-2}\leqslant r\leqslant 1\}$ is a one-time shift of the flow $\xi_A^\tau(w,z)=\left(\lambda^{\tau}w,2^{-\tau}z\right)$ and $L_0$ is transversal to trajectories of this flow. Then $G_0\cap\tilde f_A(G_0)=\emptyset$ and $G_0$ with $G_1=\tilde f_A(G_0)$ bounds a compact subset $K_A$ which is a fundamental domain of the  restriction of the diffeomorphism $f$ to $V_f$. 

Let us construct a diffeomorphism $\eta_A:K_A\to  S_{{A}}\times[0,1]$. For this aim let us define the flow $\zeta_A^\tau(w,z)=\left(w,4^{-\tau}z\right)$ on the set $\{(r,z): r\geqslant 1\}$. Finitely, let us $E^\tau_A$ be a confluence of the flows $\xi^\tau_A$ and $\zeta^\tau_A$ on $K_A$ (see Fig. \ref{signu}) whose trajectories transversally intersects $G_0$ and $G_1$ at exactly one point. Then  for every point $s\in G_0$ it is correctly defined a time $\tau(s)$ such that $E^{\tau(s)}(s)\in G_1$. Thus the desired diffeomorphism $\eta_A$ assigns for a point $(w,s)=E^{t\tau(s)}(s),\,s\in S_A,\,t\in[0,1]$ from $K_A$ the point $(s,t)\in S_A\times[0,1]$.

Define a diffeomorphism $\Psi_A:S_A\to S_A$ by the formula $$\eta_Af\eta_A^{-1}(s,0)=(\Psi_A(s),1).$$ By the construction $\Psi_A$ satisfies all requirements of the lemma.
\end{demo}

Using similar notation with index $R$ for the  repeller of the diffeomorphism $f$ we conclude that:
\begin{itemize}
\item $\Sigma_R$ is a compact orientable surface of genus $g_R$ with $m_R$ boundary components and negative Euler characteristic
\item $R$ has a trapping  neighborhood $U_{R}$ of the form $\Sigma_R \times [-1,1] $, where $\Sigma_R =\Sigma_R \times \{0 \} $ is a trapping neighborhood of repeller $R$ as a repeller of a 2-diffeomorphism $\psi_{_R}$ and diffeomorphism $ f^{-1}_{_R}=f^{-1}|_{U_R}:U_R\to f^{-1}(U_R)$ has a form  $f^{-1}_{_R} (w, z) = (\psi_{_R} (w), z/2): \Sigma_R \times [-1,1 ]\to f^{-1}(\Sigma_R)\times [-1/2,1/2] $;
\item $S_{R}$ is a closed surface of a positive genus $\rho_R=2g_R+m_R-1>1$ formed by gluing two copies of $\Sigma_R$ along the boundary components, $C_R$ is the coinciding copies of $\partial\Sigma_R$ in $S_R$;
\item diffeomorphism $f|_{V_f}$ is smoothly conjugate with a diffeomorphism $\chi_R:S_{R}\times\mathbb R\to S_{R}\times\mathbb R$ given by the formula 
$$\chi_R(s,r)=(\Psi_R(s),r+1),$$ where $\Psi_R:S_R\to S_R$ coincides with $\psi_{_R}$ on each copy of $\Sigma_R$ outside a neighborhood of $\partial\Sigma_R$ and is a contraction to $\partial\Sigma_R$ in the neighborhood being identical on $\partial\Sigma_R$.
\end{itemize} 

Let us introduce  orbit spaces $\hat V_f=V_f/f$,  $\hat V_A=(S_{A}\times\mathbb R)/\chi_A$, $\hat V_R=(S_{R}\times\mathbb R)/\chi_R$ and the natural projections $p_f:V_f\to\hat V_f$, $p_A:S_{A}\times\mathbb R\to\hat V_A$, $p_R:S_{R}\times\mathbb R\to\hat V_R$. Also let $T_A=C_A\times\mathbb R,\,\hat T_A=p_A(T_A)$,  $T_R=C_R\times\mathbb R,\,\hat T_R=p_R(T_R)$. 

\begin{lemm}\label{A=R} For every $\Omega$-stable diffeomorphism $f:M^3\to M^3$  with one-dimensional canonically embedded surface attractor-repeller the following is true:
\begin{enumerate}
\item $\rho_A=\rho_R$;
\item $\hat T_A$ $(\hat T_R)$ is a collection of $m_A$ $(m_R)$ incompressible tori which form  characteristic submanifold in $\hat V_A$ $(\hat V_R)$, that is $\hat V_A$ $(\hat V_R)$ should be cut along these tori to yield pieces that each have geometric structures (the JSJ-decomposition);
\item $m_A=m_R,\,g_A=g_R$.
\end{enumerate} 
\end{lemm}
\begin{demo} $ $

1. It follows from Lemma \ref{sia} that $V_f$ is diffeomorphic to $ S_{A}\times\mathbb R$ and $S_{R}\times\mathbb R$ simultaneously, then $ S_{A}\times\mathbb R$ and $S_{R}\times\mathbb R$ are diffeomorphic. Hence, the fundamental group $\pi_1(S_{A}\times\mathbb R)$ and $\pi_1(S_{R}\times\mathbb R)$ are isomorphic, that implies $\rho_A=\rho_R$.

2. By the construction $\Psi_A$ ($\Psi_R$) is a  diffeomorphism reducible by curves $C_A$ ($C_R$). As $\Psi_A$ ($\Psi_R$) is identical on $C_A$ ($C_R$) then  $\hat T_A$ ($\hat T_R$) are $m_A$ ($m_R$) pairwice disjoint tori. 
It follows from the Nilsen-Thurston theory \cite{Thu} that are incompressible and form  characteristic submanifold in $\hat V_A$ ($\hat V_R$), that is $\hat V_A$ ($\hat V_R$) should be cut along these tori to yield pieces that each have geometric structures (the JSJ-decomposition). 

3. It follows from Lemma \ref{sia} that $\hat V_A$ and $\hat V_R$ are diffeomorphic. As  
characteristic submanifold is unique up to isotopy (see, for example \cite{Boh}) then $m_A=m_R$.
\end{demo}

Thus, without loss of generality, we can assume that $\Sigma_A=\Sigma_R=\Sigma$, $S_A=S_R=S$ and $C_A=C_R=C$. Let $V=S\times \mathbb R$. It follows from Lemma \ref{sia} that there are diffeomorphisms $\eta_A,\,\eta_R:V_f\to V$ such that $$\eta_Af=\chi_A\eta_A,\,\eta_Rf=\chi_R\eta_R.$$ Then the diffeomorphism 
$\eta=\eta_A\eta_R^{-1}:V\to V$ possesses the property $$\eta\chi_R=\chi_A\eta.$$

Thus the set $V$ is foliated by  foliations $F^s=\{\eta_A(W^s_a),\,a\in A\}$ and $F^u=\{\eta_R(W^u_r),\,r\in R\}$. To prove  Theorem \ref{the-hyp} we have to show that the foliations $F^s$ and $\eta(F^u)$ are not transversal. We are going to do it using the universal covering space for $V$.

\subsection{Lifting wandering dynamics to the universal covering space}

Let $\mathbb G_{S}$ be a group of the hyperbolic isometrics of the Lobachevsky plane $\mathbb U$,  which is isomorphic to the fundamental group of the surface  $S$ (see Section \ref{desc}) and $S=\mathbb U/\mathbb G_S$. Then the group $\mathbb G_{V}=\{g\times id:\,g\in\mathbb G_{S}\}$ is a group of  isometrics of $\mathbb U\times\mathbb R$,  which is isomorphic to the fundamental group of the surface  $V$  and $V=(\mathbb U\times\mathbb R)/\mathbb G_{V}$. Denote by $q:\mathbb U\times\mathbb R\to V$ the natural projection and by $\tilde A$ a connected component of the set $q^{-1}(A)$ for a connected subset $A\subset V$. Let  $H:V\to V$ be a homeomorphism and $\tilde H:\mathbb U\times\mathbb R\to\mathbb U\times\mathbb R$ be its a lift. Then $\tilde H$ induces  the automorphism  $\tau_{\tilde H}$ of the group $\mathbb G_{S}$ which assigns the element $g'=\tau_{\tilde H}(g)$ to an element $g$ by $$(g'\times id)\circ\tilde H=\tilde H\circ(g\times id).$$
By the construction every lifts $\tilde\chi_A,\,\tilde\chi_R:\mathbb U\times\mathbb R\to\mathbb U\times\mathbb R$ of the diffeomorphisms $\chi_A,\chi_R:V\to V$ have forms  $\tilde\chi_A(\bar s,r)=(\bar\Psi_A(\bar s),r+1),\tilde\chi_R(\bar s,r)=(\bar\Psi_R(\bar s),r+1),\,(\bar s,r)\in V$, where $\bar\Psi_A,\bar\Psi_R:\mathbb U\to\mathbb U$ are lifts of $\Psi_A,\Psi_R:S\to S$. Herewith, $\tau_{\tilde \chi_A}=\tau_{\bar\Psi_A},\,\tau_{\tilde \chi_R}=\tau_{\bar\Psi_R}$. 

\begin{lemm}\label{pi1} The diffeomorphisms $\Psi_A$ and $\Psi_R$ are $\pi_1$-conjugate.\footnote{Notice that this fact follows independently from the other side. The manifold $\hat V_A$ ($\hat V_R$) is the mapping torus  $S_{A}\times[0,1]/_\sim$ ($S_{R}\times[0,1]/_\sim$) with the monodromy $\Psi_A$ ($\Psi_R$), that is $(s,1)\sim(\Psi_A(s),0)$ ($(s,1)\sim(\Psi_R(s),0)$). As the  mapping tori $\hat V_A$ ($\hat V_R$) are homeomorphic then the  monodromies $\Psi_A$ ($\Psi_R$) are $\pi_1$-conjugate (see, for example, \cite{Kuz}).}
\end{lemm}
\begin{demo} Recall that $\eta\chi_R=\chi_A\eta.$ Let $\tilde\chi_A:V\to V$ be a lift of $\chi_A$ and $\tilde\eta:V\to V$ be a lift of $\eta$. Then $\tilde\chi_R=\tilde\eta^{-1}\tilde\chi_A\tilde\eta:V\to V$ be a lift of $\chi_R$. Thus $\tau_{\tilde\eta}\tau_{\tilde\chi_R}=\tau_{\tilde\chi_A}\tau_{\tilde\eta}$ and, hence, $$\tau_{\tilde\eta}\tau_{\bar\Psi_R}=\tau_{\bar\Psi_A}\tau_{\tilde\eta}.$$
\end{demo}

\begin{figure}[h]\centerline{\includegraphics [width = 13 true cm]{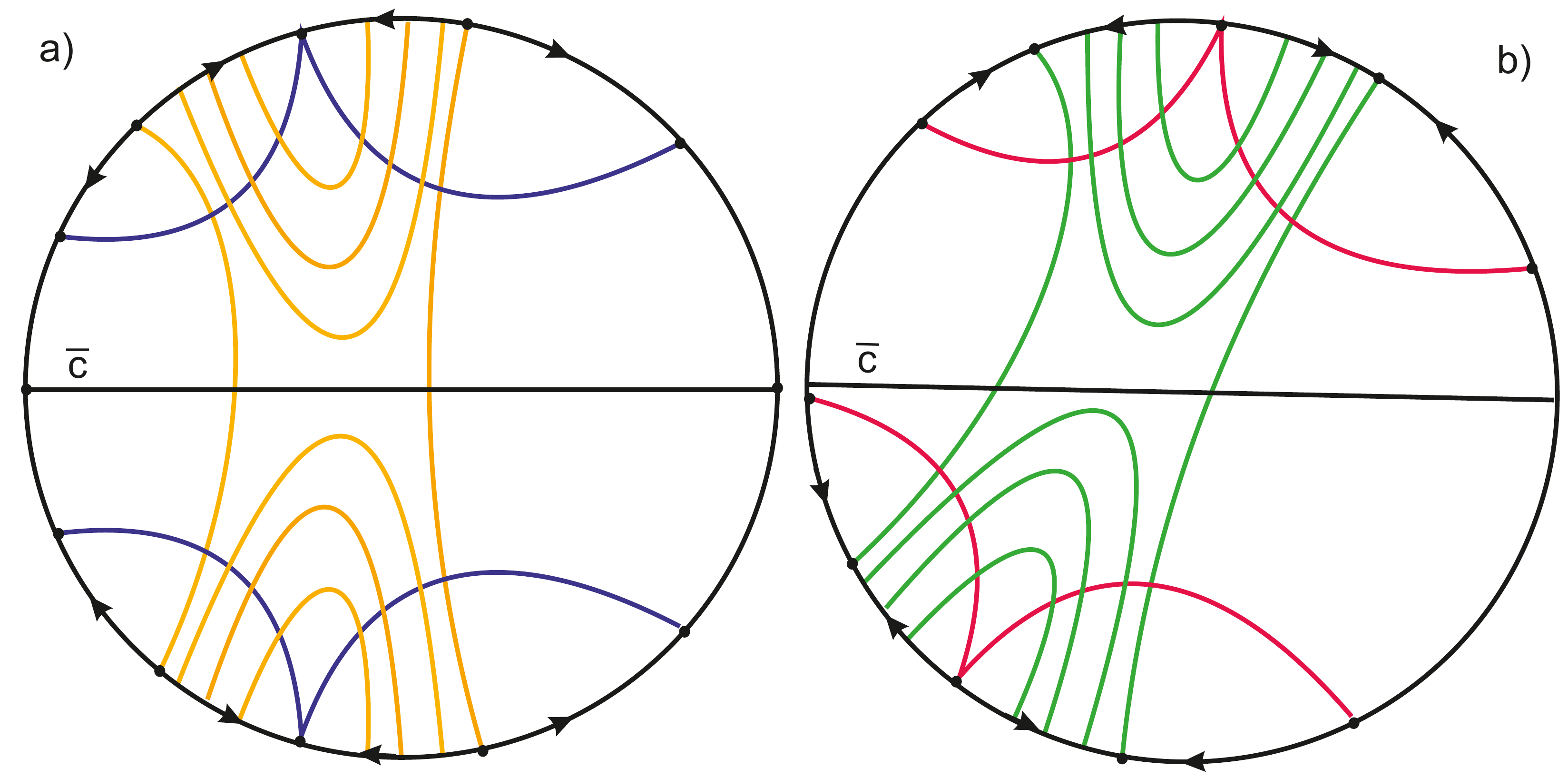}}\caption{Dynamics of the maps: a) $\bar\Psi_A$ on $\mathbb U$ and $\tau^e_{\bar\Psi_A}$ on $\mathbb E$; b) $\bar\Psi_R$ on $\mathbb U$ and $\tau^e_{\bar\Psi_R}$ on $\mathbb E$}\label{RWA}\end{figure} 

Let us fixed a curve $c$ in $C$, a lift $\bar\Psi_A:\mathbb U\to\mathbb U$ such that $\bar\Psi_A(\bar s)=\bar s$ for every point $\bar s\in\bar c$ (see Fig. \ref{RWA} a)) and the lift $\tilde\chi_A(\bar s,r)=(\bar\Psi_A(\bar s),r+1)$. Due to Lemma \ref{pi1}, without loss of generality we can assume that $\tau_{\tilde\eta}=id$ in the equality\footnote{For this aim it is enough to consider a diffeomorphism $g\Psi_Rg^{-1}$ instead $\Psi_R$ where $g:S\to S$ has a lift $\bar g:\mathbb U\to\mathbb U$ such that $\tau_{\bar g}=\tau_{\tilde\eta}$} 
$\tau_{\tilde\eta}\tau_{\tilde\chi_R}=\tau_{\tilde\chi_A}\tau_{\tilde\eta}$. Then the diffeomorphism $\bar\Psi_R$ in the lift $\tilde\chi_R(\bar s,r)=(\bar\Psi_R(\bar s),r+1)$ has a dynamics represented on Figure \ref{RWA} b). 

\section{There are no structural stable 3-diffeomorphisms with dynamic one-dimensional canonically embedded surface attractor-repeller}
In this section we prove  Theorem \ref{the-hyp}.

\begin{demo} Let $\tilde F^s=q^{-1}(F^s),\,\tilde F^u=q^{-1}(F^u)$ (see Fig.  \ref{Fs}). Due to results from the previous section it is enough to show that the foliations $\tilde F^s,\,\tilde\eta(\tilde F^u)$ are not transversal.
\begin{figure}[h]\centerline{\includegraphics [width = 10 true cm]{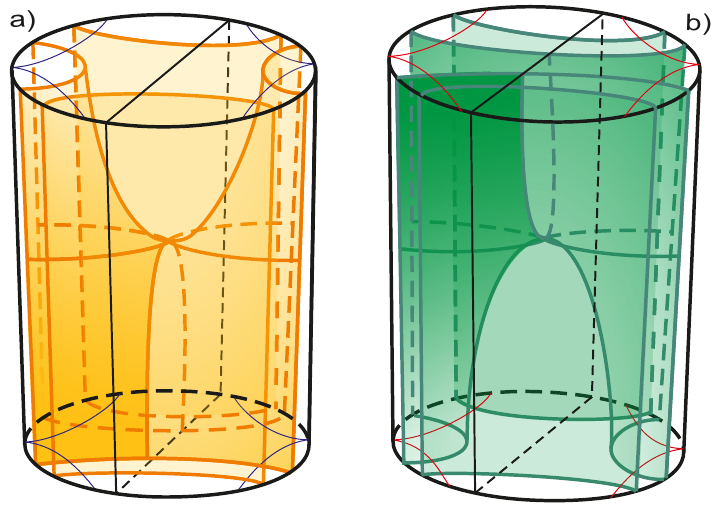}}\caption{Leaves of the foliations on $V$: a) $\tilde F^s$; b) $\tilde F^u$}\label{Fs}\end{figure}

As $\tau_{\tilde\eta}=id$ then $\tilde\eta$ preserves asymptotic behaviour of leaves of the foliation $\tilde F^u$ and, hence, mutual configuration of the leaves $\tilde w^s,\,\tilde w^u$ of the foliations $\tilde F^s,\,\tilde\eta(\tilde F^u)$ has a form represented on Figure \ref{spinst}. 
\begin{figure}[h]\centerline{\includegraphics [width = 6.5 true cm]{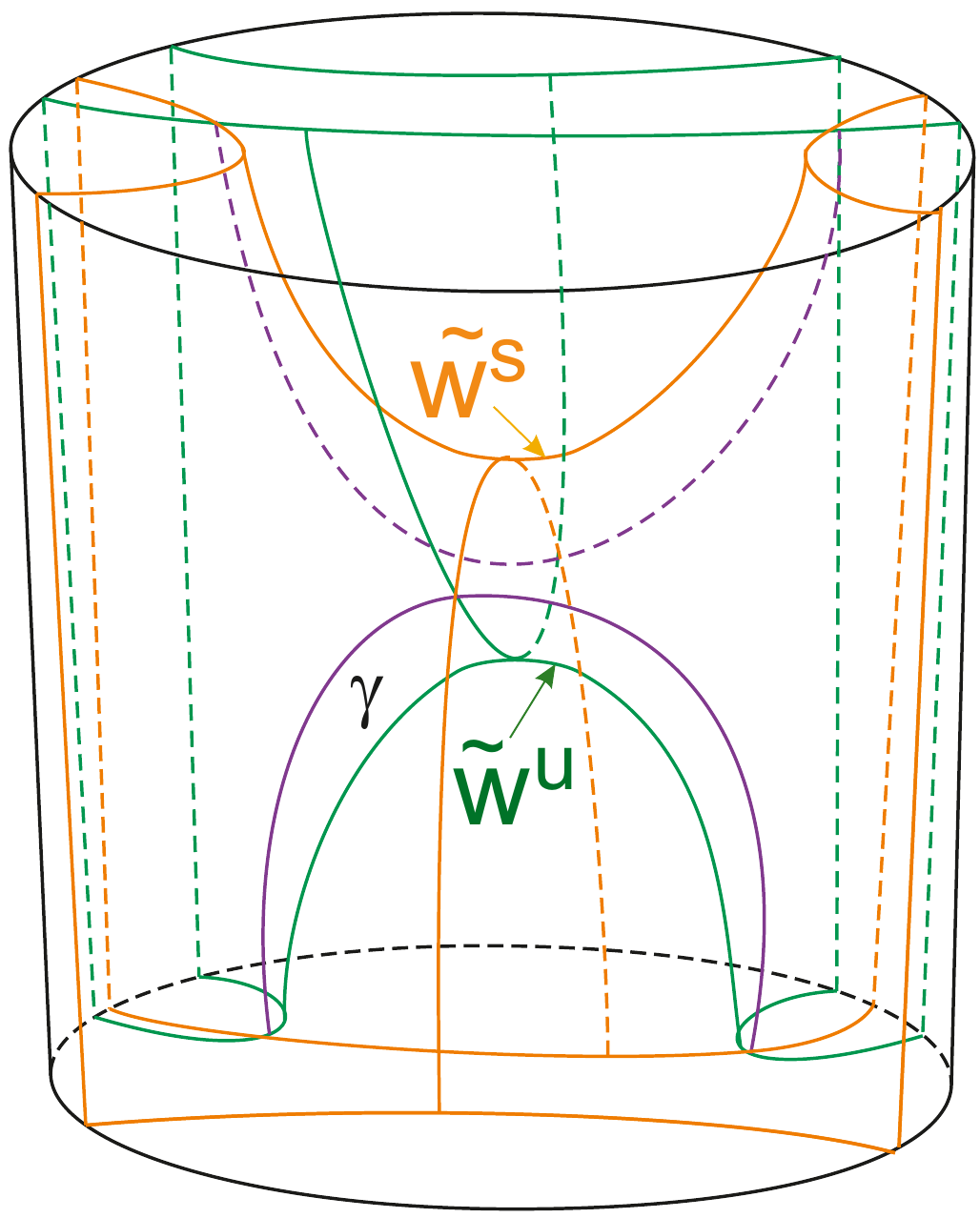}}\caption{configuration of the leaves of the foliations $\tilde F^s,\,\tilde\eta(\tilde F^u)$}\label{spinst}\end{figure} 

Every leaf $\tilde w^{s}$ is himeomorphic to 2-disc and divides $V$ into two connected components, both of them completely foliated by leaves of the foliation $\tilde  F^s$. By assumption the leaves $\tilde w^{s}$ and $\tilde w^{u}$ are transversally intersected, hence, there is a connected component $\gamma$ of the intersection $\tilde w^{u}\cap (V\setminus\tilde w^s)$ whose boundary belongs to $\tilde w^{s}$. So  there is a point  $z\in\gamma$ where $\tilde w^{u}$ has a contact with a leaf of the foliation $\tilde F^s$ that contradict to the assumption. 
\end{demo}

{\bf Acknowledgments.} This work was financially supported by the Russian Science Foundation (project 21-11-00010), except for  section \ref{dynA-R} devoted to investigation of dynamics one-dimensional canonically embedded surface attractor-repeller for 3-diffeomorphism, written with financial support from the HSE Laboratory of Dynamic Systems and Applications, the Ministry of Science and Higher Education of the Russian Federation, Contract No. 075-15-2019-1931.


\begin{thebibliography}{99}

\bibitem{Ba} A. Banyaga. The structure of the group of equivariant diffeomorphism. Topology, 16 (1977), 279--283.
\bibitem{BGPY} 
Barinova M., Grines V., Pochinka O., Yu B. Existence of an energy function for three-dimensional chaotic ``sink-source'' cascades // Chaos. 2021. Vol. 31. No. 6.
\bibitem{Boh} Francis Bonahon. Geometric structures on 3-manifolds. In Handbook of geometric topology, pages 93--164.
North-Holland, Amsterdam, 2002.
\bibitem{Bonatti2010}  Bonatti Ch., Guelman N. Axiom A diffeomorphisms derived from Anosov flows. J. Mod. Dyn. 4 (2010), no. 1, 1--63.
\bibitem{Brown2010} Brown A. Nonexpanding attractors: conjugacy to algebraic models and classification in 3-manifolds. // Journal of Modern Dynamics. 2010. V. 4.  517--548.
\bibitem{Gib72}  Gibbons J.C. One-Dimensional basic sets in the three-sphere. Trans. of the
Amer. Math. Soc.  Volume 164, February 1972, 163-178.

\bibitem{Grines1975} V. Grines. On topological conjugacy of diffeomorphisms of a
two-dimensional manifold onto one-dimensional orientable basic
sets I, Transactions of the Moscow Mathematical Society 32,
31-56 (1975).
\bibitem{GrMi} V. Z. Grines. On the topological classification of structurally stable diffeomorphisms of surfaces with one-dimensional attractors and repellers. Sb. Math., 188:4 (1997), 537--569.
\bibitem{GrLePo} {Grines V., Levchenko Y., Medvedev V., Pochinka O.} (2015). The topological classification of structural stable 3-diffeomorphisms with two-dimensional basic sets, (Nonlinearity) 28, No.11, 4081--4102.
\bibitem{GrMeZh} V.Z.~Grines, V.S. Medvedev,
E.V. Zhuzhoma.  On surface attractors and repellers on 3-manifolds. Mat. zam. 2005. V. 78 (6). 821--834.
\bibitem{GrinesZh2005} Grines V., Zhuzhoma E. On structurally stable diffeomorphisms with codimension one expanding attractors. // Trans. Amer. Math. Soc. 2005. V. 357 (2). 617--667.
\bibitem{GrPo} V. Grines, T. Medvedev, and O. Pochinka, Dynamical Systems
on 2- and 3-Manifolds (Switzerland: Springer, 2016).
\bibitem{HTh} Handel, M.; Thurston, W. P. (1985). New proofs of some results of Nielsen. Advances in Mathematics. 56 (2): 173--191.
\bibitem{JNW} B. Jiang, Y. Ni, S. Wang, 3-manifolds that admit knotted solenoids as attractors, Trans.
Amer. Math. Soc. 356 (2004), no. 11, 4371--4382.
\bibitem{Kuz} N. F. Kuzennyi. Isomorphism of semidirect products. Ukrainian Mathematical Journal, 26 (5): 543--547, 1974.
\bibitem{MeZhu2005} Medvedev V., Zhuzhoma E. Journal of Dynamical and Control Systems. 2005. V. 11 (3). 405--411. 
\bibitem{Plykin1971} R. Plykin, The topology of basic sets for smale diffeomorphisms, Math. USSR-Sb. 13, 301--312 (1971).
\bibitem{PlZh}  Plykin R. V.,  Zhirov A. Yu. On the relationship between one-dimensional hyperbolic attractors of surface diffeomorphisms and generalized pseudo-Anosov diffeomorphisms.  Mat. notes. 1995.
V. 58. № 1. P. 149--152.
\bibitem{RoWi} {Robinson R.C., Williams R.F.} (1973). Finite Stability is not generic (Dynamical Systems ((Proc. Sympos., Univ. Bahia, Salvador, 1971)) Academic Press, New York), 451--462.
\bibitem{Shi2014} Shi Yi. Partially hyperbolic diffeomorphisms on Heisenberg nilmanifolds and
holonomy maps. C. R. Math. Acad. Sci. Paris. 2014. V. 352, No.9, 743--747.
\bibitem{Thu} W. Thurston, On the geometry and dynamics of diffeomorphisms of surfaces, Bulletin of the American Mathematical Society, vol. 19 (1988), pp. 417--431.
\bibitem{Wi74} R. Williams. Expanding attractors. Inst. Hautes Etudes Sci. Publ. Math. 1974. V. 43. 169--203. 


\end{thebibliography}
\end{document}